\documentclass[12pt]{article}
\usepackage{amsmath,amsfonts,amssymb,amsthm}
\input amssym.def
\begin{document}
\def \Z{\Bbb Z}
\def \C{\Bbb C}
\def \R{\Bbb R}
\def \Q{\Bbb Q}
\def \N{\Bbb N}

\def \A{{\mathcal{A}}}
\def \D{{\mathcal{D}}}
\def \E{{\mathcal{E}}}
\def \E{{\mathcal{E}}}
\def \H{\mathcal{H}}
\def \S{{\mathcal{S}}}
\def \wt{{\rm wt}}
\def \tr{{\rm tr}}
\def \span{{\rm span}}
\def \Res{{\rm Res}}
\def \Der{{\rm Der}}
\def \End{{\rm End}}
\def \Ind {{\rm Ind}}
\def \Irr {{\rm Irr}}
\def \Aut{{\rm Aut}}
\def \GL{{\rm GL}}
\def \Hom{{\rm Hom}}
\def \mod{{\rm mod}}
\def \ann{{\rm Ann}}
\def \ad{{\rm ad}}
\def \rank{{\rm rank}\;}
\def \<{\langle}
\def \>{\rangle}

\def \g{{\frak{g}}}
\def \h{{\hbar}}
\def \k{{\frak{k}}}
\def \sl{{\frak{sl}}}
\def \gl{{\frak{gl}}}

\def \be{\begin{equation}\label}
\def \ee{\end{equation}}
\def \bex{\begin{example}\label}
\def \eex{\end{example}}
\def \bl{\begin{lem}\label}
\def \el{\end{lem}}
\def \bt{\begin{thm}\label}
\def \et{\end{thm}}
\def \bp{\begin{prop}\label}
\def \ep{\end{prop}}
\def \br{\begin{rem}\label}
\def \er{\end{rem}}
\def \bc{\begin{coro}\label}
\def \ec{\end{coro}}
\def \bd{\begin{de}\label}
\def \ed{\end{de}}

\newcommand{\m}{\bf m}
\newcommand{\n}{\bf n}
\newcommand{\nno}{\nonumber}
\newcommand{\nord}{\mbox{\scriptsize ${\circ\atop\circ}$}}
\newtheorem{thm}{Theorem}[section]
\newtheorem{prop}[thm]{Proposition}
\newtheorem{coro}[thm]{Corollary}
\newtheorem{conj}[thm]{Conjecture}
\newtheorem{example}[thm]{Example}
\newtheorem{lem}[thm]{Lemma}
\newtheorem{rem}[thm]{Remark}
\newtheorem{de}[thm]{Definition}
\newtheorem{hy}[thm]{Hypothesis}
\makeatletter \@addtoreset{equation}{section}
\def\theequation{\thesection.\arabic{equation}}
\makeatother \makeatletter

\begin{center}
{\Large \bf Vertex $F$-algebras and their $\phi$-coordinated
modules}
\end{center}

\begin{center}
{Haisheng Li\footnote{Partially supported by NSF grant DMS-0600189}\\
Department of Mathematical Sciences\\
Rutgers University, Camden, NJ 08102}
\end{center}

\begin{abstract}
In this paper, for every one-dimensional formal group $F$ we
formulate and study a notion of vertex $F$-algebra and a notion of
$\phi$-coordinated module for a vertex $F$-algebra where $\phi$ is
what we call an associate of $F$. In the case that $F$ is the
additive formal group, vertex $F$-algebras are exactly ordinary
vertex algebras.  We give a canonical isomorphism between the
category of vertex $F$-algebras and the category of ordinary vertex
algebras. Meanwhile, for every formal group we completely determine
its associates. We also study $\phi$-coordinated modules for a
general vertex $\Z$-graded algebra $V$ with $\phi$ specialized to a
particular associate of the additive formal group and we give a
canonical connection between $V$-modules and $\phi$-coordinate
modules for a vertex algebra obtained from $V$ by Zhu's
change-of-variables theorem.
\end{abstract}

\section{Introduction}
In a series of papers (see \cite{li-qva1}, \cite{li-qva2},
\cite{li-infinity}, \cite{li-hqva}, \cite{li-phi-module}), we have
been extensively investigating various vertex algebra-like
structures naturally arising from Yangians and quantum affine
algebras. One of the main goals is to solve the problem (see
\cite{fj}, \cite{efk}), to associate ``quantum vertex algebras" to
quantum affine algebras. Partly motivated by Etingof-Kazhdan's
theory of quantum vertex operator algebras (see \cite{ek}), in
\cite{li-qva1} we developed a theory of (weak) quantum vertex
algebras and we established a conceptual construction of weak
quantum vertex algebras and their modules. To associate weak quantum
vertex algebras to quantum affine algebras, in \cite{li-phi-module}
we furthermore developed a theory of what we called
$\phi$-coordinate quasi modules for weak quantum vertex algebras and
established a general construction, where $\phi$ is what we called
an associate of the one-dimensional additive formal group.

In this paper, we continue to formulate and study notions of
nonlocal vertex $F$-algebra and vertex $F$-algebra with $F$ a
one-dimensional formal group. It is shown that the category of
vertex $F$-algebras is canonically isomorphic to the category of
ordinary vertex algebras. We also study a notion of
$\phi$-coordinated quasi module for a vertex $F$-algebra $V$ with
$\phi$ an associate of $F$ in the sense of \cite{li-phi-module}, and
we show that the category of $\phi$-coordinated quasi $V$-modules is
canonically isomorphic to the category of $\hat{\phi}$-coordinated
quasi $\hat{V}$-modules, where $\hat{V}$ is a certain ordinary
vertex algebra and $\hat{\phi}$ is a certain associate of the
additive formal group.

In the following, we give a more detailed account of this paper. Let
us start with the notion of vertex algebra, using one of several
equivalent definitions (cf. \cite{ll}). A {\em vertex algebra}  is a
vector space $V$, equipped with a linear map
$$Y(\cdot,x): V\rightarrow \Hom (V,V((x)))\subset (\End
V)[[x,x^{-1}]]$$ and equipped with a distinguished vector ${\bf
1}\in V$, satisfying the following conditions:\\
 (I) {\em Vacuum and creation property}
$$Y({\bf 1},x)v=v,\ \ Y(v,x){\bf 1}\in V[[x]]\ \mbox{ and }\
\lim_{x\rightarrow 0}Y(v,x){\bf 1}=v\ \ \mbox{ for }v\in V.$$ (II)
{\em Weak commutativity:} For $u,v\in V$, there exists $k\in \N$
such that
\begin{eqnarray}\label{weakcomm}
(x_{1}-x_{2})^{k}Y(u,x_{1})Y(v,x_{2})=(x_{1}-x_{2})^{k}Y(v,x_{2})Y(u,x_{1}).
\end{eqnarray}
(III) {\em Weak associativity:} For $u,v,w\in V$, there exists $l\in
\N$ such that
\begin{eqnarray}\label{weakassoc}
(x_{0}+x_{2})^{l}Y(u,x_{0}+x_{2})Y(v,x_{2})w =(x_{0}+x_{2})^{l}
Y(Y(u,x_{0})v,x_{2})w.
\end{eqnarray}
Note that this weak associativity relation can be viewed as an
algebraic version of what physicists often call operator product
expansion, while the weak commutativity relation is a version of
what physicists call locality. Taking out the weak commutativity
axiom from the above definition, one arrives at the notion of what
was called nonlocal vertex algebra in \cite{li-qva1}. (Nonlocal
vertex algebras are the same as weak $G_{1}$-vertex algebras in
\cite{li-g1} and are also essentially field algebras in \cite{bk}.)
While vertex algebras are analogs of commutative associative
algebras, nonlocal vertex algebras are analogs of noncommutative
associative algebras.

One-dimensional formal groups (over $\C$) (cf. \cite{haz}) are
formal power series $F(x,y)\in \C[[x,y]]$, satisfying
$$F(x,0)=x,\ \ \ \ F(0,y)=y,\ \ \ \ F(x,F(y,z))=F(F(x,y),z),$$
among which the simplest is the {\em additive formal group} $F_{\rm
a}(x,y)=x+y$. What we called associates of a formal group $F$ are
formal series $\phi(x,z)\in \C((x))[[z]]$, satisfying
\begin{eqnarray}
\phi(x,0)=x,\ \ \ \ \
\phi(\phi(x,z_{1}),z_{2})=\phi(x,F(z_{1},z_{2})).
\end{eqnarray}
(To a certain extent, an associate of $F$ to a formal group $F$  is
like a $G$-set to a group $G$.) Associates of $F_{\rm a}$ have been
completely determined in \cite{li-phi-module}, where two particular
examples are $\phi(x,z)=x+z=F_{\rm a}(x,z)$ and $\phi(x,z)=xe^{z}$.
In this present paper, we furthermore have completely determined the
associates of a general one-dimensional formal group $F$.

Let $F$ be a one-dimensional formal group. We define a notion of
(nonlocal) vertex $F$-algebra by simply replacing $x_{0}+x_{2}\
(=F_{\rm a}(x_{0},x_{2}))$ in the weak associativity relation
(\ref{weakassoc}) with $F(x_{0},x_{2})$, i.e., for $u,v,w\in V$,
there exists $l\in \N$ such that
\begin{eqnarray}\label{Fweakassoc}
F(x_{0},x_{2})^{l}Y(u,F(x_{0},x_{2}))Y(v,x_{2})w =F(x_{0},x_{2})^{l}
Y(Y(u,x_{0})v,x_{2})w.
\end{eqnarray}
Recall that the {\em logarithm} of the formal group $F$ (cf.
\cite{haz}) is the unique series $f(x)\in x\C[[x]]$ such that
$f'(0)=1$ and $f(F(x,y))=f(x)+f(y).$ Let $(V,Y,{\bf 1})$ be a (resp.
nonlocal) vertex algebra. For $v\in V$, set $Y_{F}(v,x)=Y(v,f(x))$.
We show that $(V,Y_{F},{\bf 1})$ is a (resp. nonlocal) vertex
$F$-algebra. Furthermore, we show that this gives rise to a
canonical isomorphism between the category of (resp. nonlocal)
vertex algebras and the category of (resp. nonlocal) vertex
$F$-algebras.

Given a vertex $F$-algebra $V$, we define a notion of $V$-module in
the obvious way. More generally we define a notion of
$\phi$-coordinated $V$-module with $\phi$ an associate of $F$, where
the defining weak associativity relation reads as
\begin{eqnarray}\label{Fweakassoc-module}
\left((x_{1}-x_{2})^{k}Y(u,x_{1})Y(v,x_{2})\right)|_{x_{1}=\phi(x_{2},x_{0})}
=(\phi(x_{2},x_{0})-x_{2})^{k} Y(Y(u,x_{0})v,x_{2}),
\end{eqnarray}
where $k$ is a certain nonnegative integer. Furthermore, we study a
general construction of nonlocal vertex $F$-algebras and their
$\phi$-coordinated modules. Let $W$ be a general vector space and
set $\E(W)=\Hom (W,W((x)))$. Given an associate $\phi$ of $F$, for
any pair $(a(x),b(x))$ in $\E(W)$, satisfying a certain condition,
we define $a(x)_{n}^{\phi}b(x)\in \E(W)$ for $n\in \Z$ in terms of
the generating function
$$Y_{\E}(a(x),z)b(x)=\sum_{n\in \Z}a(x)_{n}^{\phi}b(x)z^{-n-1}$$
roughly by
\begin{eqnarray}
``Y_{\E}^{\phi}(a(x),z)b(x)=\left(a(x_{1})b(x)\right)|_{x_{1}=\phi(x,z)}."
\end{eqnarray}
(See Section 4 for the precise definition.) It is proved that every
what we call compatible subset of $\E(W)$ generates a nonlocal
vertex $F$-algebra with $W$ as a canonical $\phi$-coordinated
module. This generalizes the corresponding result of \cite{li-qva1}.

Interestingly, the notion of $\phi$-coordinated module for vertex
operator algebras has an intrinsic connection with Zhu's work on
change-of-variables (also see the next paragraph). More
specifically, let $(V,Y,{\bf 1},\omega)$ be a vertex operator
algebra in the sense of \cite{flm}. For $v\in V$, set
$Y[v,x]=Y(e^{xL(0)}v,e^{x}-1)$, where $Y(\omega,x)=\sum_{n\in
\Z}L(n)x^{-n-2}$. Set $\tilde{\omega}=\omega-\frac{c}{24}{\bf 1}$,
where $c$ is the central charge of $V$. It was proved by Zhu
(\cite{zhu1}, \cite{zhu2}; cf. \cite{lep1}) that $(V,Y[\cdot,x],{\bf
1},\tilde{\omega})$ carries the structure of a vertex operator
algebra. Now, let $(W,Y_{W})$ be a $V$-module. For $v\in V$, set
$$X_{W}(v,z)=Y_{W}(z^{L(0)}v,z)\in (\End W)[[z,z^{-1}]].$$
 We show that $(W,X_{W})$ carries
the structure of a $\phi$-coordinated module with $\phi(x,z)=xe^{z}$
for the vertex operator algebra $(V,Y[\cdot,x],{\bf
1},\tilde{\omega})$. In fact, what we have done in this paper is
more general with $V$ a nonlocal vertex $\Z$-graded algebra.

In a series of papers (see \cite{lep2}, \cite{lep3}, \cite{lep4}),
Lepowsky has extensively studied the homogenized vertex operators
$Y(x^{L(0)}v,x)$ (with $v\in V$) for a general vertex operator
algebra $V$, and obtained a very interesting new type of Jacobi
identity, incorporating values of the Riemann zeta function at
negative integers. The substitution $x_{1}=x_{2}e^{z}$ and Zhu's
work on change-of-variables have already entered into his study. The
formal (unrigorous) relation (4.20) in \cite{lep2} is inspirational
for our introduction of the notion of $\phi$-coordinated module in
\cite{li-phi-module} and it is the main motivation of Proposition
\ref{pordinary-exp-module} in the present paper.

In \cite{b-Gva}, Borcherds studied a notion of vertex $G$-algebra,
which generalizes the notion of vertex algebra significantly. It was
pointed out therein that the ordinary vertex algebra theory was
based on the one-dimensional additive formal group and there was a
vertex algebra theory associated to every formal group. Presumably,
(nonlocal) vertex $F$-algebras belong to the family of Borcherds'
vertex $G$-algebras.

There is a very interesting paper \cite{lod}, in which to any
one-dimensional formal group, Loday associated a type of algebras,
which is an operad. It was shown that for the additive formal group,
the corresponding operad is that of associative algebras with
derivation. It would be very interesting to relate \cite{lod} with
the present paper.

This paper is organized as follows: In Section 2, we study
associates of one-dimensional formal groups. In Section 3, we study
nonlocal vertex $F$-algebras and vertex $F$-algebras. In Section 4,
we study $\phi$-coordinated modules for nonlocal vertex
$F$-algebras. In Section 5, we study $\phi$-coordinated modules for
nonlocal vertex algebras with $\phi$ specialized to
$\phi(x,z)=xe^{z}$.

\section{Associates of one-dimensional formal groups}
In this section we first review the basics of one-dimensional formal
groups and we then study what we called in \cite{li-phi-module}
associates of a formal group. For every one-dimensional formal
group, we completely determine its associates.

In this paper, in addition to the standard notations $\Z$, $\Q$,
$\R$, and $\C$ for the integers, the rational numbers, the real
numbers, and the complex numbers, respectively, we use $\N$ for the
nonnegative integers. All vector spaces are assumed to be over $\C$
unless a different scalar field is specified otherwise.

We now follow \cite{haz} to recall the basics of one-dimensional
formal groups.

\bd{dformal-group} {\em Let $R$ be a commutative and associative
ring with identity. A {\em one-dimensional formal group} over $R$ is
a formal power series $F(x,y)\in R[[x,y]]$, satisfying
\begin{eqnarray}
F(x,0)=x,\ \ \ F(0,y)=y,  \ \ \ F(x,F(y,z))=F(F(x,y),z).
\end{eqnarray}} \ed

Throughout this paper, by a formal group we shall always mean a
one-dimensional formal group. A formal group $F(x,y)$ is said to be
{\em commutative} if $F(x,y)=F(y,x)$. A fact (see \cite{haz}) is
that formal groups over $R$ are all commutative unless $R$ contains
a nonzero element $a$ such that $a^{n}=0$ and $na=0$ for some
positive integer $n$.

Particular examples are the {\em additive formal group}
$$F_{\rm a}(x,y)=x+y$$
and the {\em multiplicative formal group}
$$F_{\rm m}(x,y)=x+y+xy.$$

\br{rgroup-power} {\em Let $R$ be a commutative associative algebra
over $\Q$. Set
\begin{eqnarray}
G(R[[x]])=\{ f(x)\in R[[x]]\;|\; f(0)=0,\ f'(0)=1\}.
\end{eqnarray}
For $f(x),g(x)\in G(R[[x]])$, we have $f(g(x))\in G(R[[x]])$.
Furthermore, $G(R[[x]])$ is a group with respect to the composition
operation. For $f(x)\in G(R[[x]])$, denote the inverse of $f(x)$ by
$f^{-1}(x)$. } \er

\bp{pequivalence} Let $R$ be a commutative associative algebra over
$\Q$. Let $f(x)\in G(R[[x]])$. Define
$$F(x,y)=f^{-1}(f(x)+f(y))\in R[[x,y]].$$
Then $F(x,y)$ is a formal group over $R$. Furthermore, this
association gives a bijection between $G(R[[x]])$ and the set of
formal groups over $R$. \ep

\bd{dlog} {\em Assume that $R$ is a commutative associative algebra
over $\Q$. Let $F(x,y)$ be a formal group over $R$. The {\em
logarithm of $F$}, denoted by $\log F$, is defined to be the unique
power series $f(x)\in R[[x]]$ with $f(0)=0$ and $f'(0)=1$ such that
\begin{eqnarray}\label{elog-def}
f(F(x,y))=f(x)+f(y).
\end{eqnarray}} \ed

For the additive formal group and the multiplicative formal group we
have
\begin{eqnarray}
\log F_{\rm a}=x,\ \ \ \ \log F_{\rm m}=\log (1+x),
\end{eqnarray}
where by definition
\begin{eqnarray}\label{elog-exp}
\log(1+x)=\sum_{n\ge 1}(-1)^{n-1}\frac{1}{n}x^{n}\in \Q[[x]].
\end{eqnarray}
All the basics mentioned above can be found from \cite{haz}.

For the rest of this paper, we restrict ourselves to formal groups
over $\C$. In this case, all (one-dimensional) formal groups are
commutative. Recall the following notion from \cite{li-phi-module}:

\bd{dassociate} {\em Let $F(x,y)$ be a formal group over $\C$. An
{\em associate} of $F(x,y)$ is a formal series $\phi(x,z)\in
\C((x))[[z]]$, satisfying
\begin{eqnarray}
\phi(x,0)=x,\ \ \ \ \phi(\phi(x,y),z) =\phi(x,F(y,z)).
\end{eqnarray}} \ed

To a certain extent, an associate of $F$ to a formal group $F$ is
like a $G$-set to a group $G$.  Note that every formal group
$F(x,y)$ is an associate of itself. On the other hand, it can be
readily seen that $\phi(x,z)=x$ is an associate of every formal
group.

For the additive formal group $F_{\rm a}$, we have the following
explicit construction of associates due to \cite{li-phi-module}:

\bp{pall-assoc} For $p(x)\in \C((x))$, set
$$\phi_{p(x)}(x,z)=e^{zp(x)\frac{d}{dx}}x
=\sum_{n\ge 0}\frac{z^{n}}{n!}\left(p(x)\frac{d}{dx}\right)^{n}x\in
\C((x))[[z]].$$ Then $\phi_{p(x)}(x,z)$ is an associate of $F_{\rm
a}$. Furthermore, every associate of $F_{\rm a}$ is of this form
with $p(x)$ uniquely determined. \ep

Using Proposition \ref{pall-assoc}, we obtain particular associates
of $F_{\rm a}$: $\phi_{p(x)}(x,z)=x$ with $p(x)=0$;
$\phi_{p(x)}(x,z)=x+z$ with $p(x)=1$; $\phi_{p(x)}(x,z)=xe^{z}$ with
$p(x)=x$; $\phi_{p(x)}(x,z)=x(1-zx)^{-1}$ with $p(x)=x^{2}$.

\br{rexistence} {\em  Let $g(x)\in G(\C[[x]])$,  i.e., $g(x)\in
x\C[[x]]$ with $g'(0)=1$.  We have
$$g(x)^{m}\in x^{m}\C[[x]]\ \ \ \mbox{ for }m\in \Z.$$
For any $h(x)=\sum_{m\ge k}\alpha_{m}x^{m}\in \C((x))$ (with $k\in
\Z,\; \alpha_{m}\in \C$), we have
$$h(g(x))=\sum_{m\ge k}\alpha_{m}g(x)^{m}\in x^{k}\C[[x]]\subset\C((x)).$$
It is clear that the map sending $h(x)\in \C((x))$ to $h(g(x))$ is
an automorphism of $\C((x))$, which preserves the subalgebra
$\C[[x]]$. Furthermore, we have an automorphism of $\C((x))[[z]]$,
sending $\psi(x,z)$ to $\psi(g(x),g(z))$ for $\psi(x,z)\in
\C((x))[[z]]$. We also have an automorphism of $\C((x))[[z]]$,
sending $\psi(x,z)$ to $\psi(x,g(z))$.} \er

\bp{passociate-relation} Let $F$ be a formal group, let $\phi$ be an
associate of $F$, and let $g(x)\in x\C[[x]]$ with $g'(0)=1$. Set
\begin{eqnarray*}
&&F_{g}(x,y)=g^{-1}(F(g(x),g(y)))\in \C[[x,y]],\\
&&\phi_{g}(x,z)=g^{-1}(\phi(g(x),g(z)))\in \C((x))[[z]].
\end{eqnarray*}
Then $F_{g}$ is a formal group and $\phi_{g}$ is an associate of
$F_{g}$. \ep

\begin{proof} Let $f=\log F$. That is, $F(x,y)=f^{-1}(f(x)+f(y))$. Thus
$$F_{g}(x,y)=g^{-1}\left(f^{-1}(f(g(x))+f(g(y)))\right).$$
In view of Proposition \ref{pequivalence}, $F_{g}$ is a formal group
with $\log F_{g}=f\circ g$.

First, note that $\phi_{g}(x,z)$ is a well defined element of
$\C((x))[[z]]$. Indeed, for any $h(x)=\sum_{n\ge 0}c_{n}x^{n}\in
\C[[x]]$, writing $\phi(x,z)=x+zA$ with $A\in \C((x))[[z]]$, we have
$$h(\phi(x,z))=\sum_{n\ge 0}c_{n}\phi(x,z)^{n}
=\sum_{j\ge 0}\left(\sum_{n\ge 0}\binom{n}{j}c_{n}x^{n-j}A^{j}
\right)z^{j}\in \C((x))[[z]].$$ Furthermore, we have
$$\phi_{g}(x,0)=g^{-1}(\phi(g(x),g(0)))=g^{-1}(\phi(g(x),0))=g^{-1}(g(x))=x$$
and
\begin{eqnarray*}
&&\phi_{g}(\phi_{g}(x,y),z))=g^{-1}\left(\phi(g(\phi_{g}(x,y)),g(z))\right)
=g^{-1}\left(\phi(\phi(g(x),g(y)),g(z))\right)\\
&=& g^{-1}\left(\phi(g(x),F(g(y),g(z)))\right)
=\phi_{g}(x,F_{g}(y,z)).
\end{eqnarray*}
This proves that $\phi_{g}(x,z)$ is an associate of $F_{g}$.
\end{proof}

As an immediate consequence of Proposition
\ref{passociate-relation}, we have (cf. Proposition
\ref{pequivalence}):

\bc{cclass-Fassociates} Let $F$ be a formal group over $\C$ with
$\log F=f$. For any associate $\phi(x,z)$ of $F_{\rm a}$, set
$$\bar{\phi}(x,z)=f^{-1}(\phi(f(x),f(z)))\in \C((x))[[z]].$$
Then $\bar{\phi}(x,z)$ is an associate of $F$. Furthermore, this
gives a 1-1 correspondence between the set of associates of $F_{\rm
a}$ and the set of associates of $F$. \ec

The following is another connection between the associates of a
general formal group and the associates of the additive formal group
$F_{\rm a}$:

\bp{psimple-connection} Let $F$ be a formal group over $\C$ and let
$g(x)\in x\C[[x]]$ with $g'(0)=1$. Then for any associate
$\phi(x,z)$ of $F$, $\phi(x,g(z))$ is an associate of $F_{g}$, where
$F_{g}$ is given as in Proposition \ref{passociate-relation}.
Furthermore, the map, sending $\phi(x,z)\in \C((x))[[z]]$ to
$\phi(x,g(z))$, gives a 1-1 correspondence between the set of
associates of $F$ and the set of associates of $F_{g}$. \ep

\begin{proof} As $\phi(x,z)\in \C((x))[[z]],\ g(x)\in x\C[[x]]$, we see that
$\phi(x,g(z))$ exists in $\C((x))[[z]]$. We have
$\phi(x,g(0))=\phi(x,0)=x$ and
$$\phi(\phi(x,g(z_{1})),g(z_{2}))=\phi(x,F(g(z_{1}),g(z_{2})))
=\phi(x,g(F_{g}(z_{1},z_{2}))).$$ Thus $\phi(x,g(z))$ is an
associate of $F_{g}$.

{}From the first part, we have that for any associate $\phi(x,z)$ of
$F$, $\phi(x,g(z))$ is an associate of $F_{g}$ and for any associate
$\psi(x,z)$ of $F_{g}$, $\psi(x,g^{-1}(z))$ is an associate of $F$.
Then the furthermore assertion follows immediately.
\end{proof}

\br{rassociate-relation} {\em Let $F$ be a formal group over $\C$
with $f=\log F$ and let $\phi$ be an associate of $F_{\rm a}$. By
Corollary \ref{cclass-Fassociates} and Proposition
\ref{psimple-connection}, both $f^{-1}(\phi(f(x),f(z)))$ and
$\phi(x,f(z))$ are associates of $F$. However, they are different in
general. To see this, let $F=F_{\rm m}$ and let $\phi(x,z)=F_{\rm
a}(x,z)=x+z$. Note that $f=\log F_{\rm m}=\log (1+x)$ and
$f^{-1}(x)=e^{x}-1$. We have
$$\phi(x,f(z))=x+f(z)=x+\log (1+z),$$ whereas
$$f^{-1}(\phi(f(x),f(z)))=f^{-1}(f(x)+f(z))=e^{\log(1+x)+\log(1+z)}-1=x+z+xz.$$
} \er

Combining Propositions \ref{psimple-connection} and \ref{pall-assoc}
we immediately have:

\bt{tF-associates} Let $F$ be a formal group over $\C$ with $\log
F=f$. Then for any $p(x)\in \C((x))$,
$\psi_{p(x)}(x,z)=e^{f(z)p(x)\frac{d}{dx}}x$ is an associate of $F$.
Furthermore, every associate of $F$ is of this form. \et

The following technical result, which generalizes a result of
\cite{li-phi-module}, plays an important role in this work:

\bl{lnonzero} Let $F$ be a formal group over $\C$ and let
$\phi(x,z)$ be an associate of $F$ with $\phi(x,z)\ne x$. Then
$q(\phi(x,z),x)\ne 0$ for any nonzero series $q(x,y)\in \C((x,y))$.
\el

\begin{proof} Let $f$ be the logarithm of $F$.
By Theorem \ref{tF-associates},
$\phi(x,z)=e^{f(z)p(x)\frac{d}{dx}}x$ for some nonzero $p(x)\in
\C((x))$. By Proposition \ref{psimple-connection},
$\phi(x,f^{-1}(z))$ is an associate of $F_{\rm a}$. As $\phi(x,z)\ne
x$, we have
 $\phi(x,f^{-1}(z))\ne x$. By Lemma 2.10 of \cite{li-phi-module}, we
 have
$q(\phi(x,f^{-1}(z)),x)\ne 0$.
 Therefore $q(\phi(x,z),x)\ne 0$.
\end{proof}

\br{rmult} {\em Recall that  $\log F_{\rm m}=\log (1+x)$. In view of
Theorem \ref{tF-associates}, the associates of $F_{\rm m}$ are of
the form
$$\phi(x,z)=e^{(\log(1+z))p(x)\frac{d}{dx}}x$$
for $p(x)\in \C((x))$.  We have particular associates of $F_{\rm
m}$: $\phi(x,z)=x$ with $p(x)=0$, $\phi(x,z)=x+\log(1+z)$ with
$p(x)=1$, $\phi(x,z)=xe^{\log (1+z)}=x(1+z)$ with $p(x)=x$,
$\phi(x,z)=x(1-x\log (1+z))^{-1}$ with $p(x)=x^{2}$. } \er

\section{Nonlocal vertex $F$-algebras and vertex $F$-algebras}
In this section we formulate and study notions of nonlocal vertex
$F$-algebra and vertex $F$-algebra with $F$ a one-dimensional formal
group. When $F=F_{\rm a}$ the additive formal group, vertex $F_{\rm
a}$-algebras are simply ordinary vertex algebras, while nonlocal
vertex $F_{\rm a}$-algebras are usual nonlocal vertex algebras. As
the main result of this section, we exhibit a canonical isomorphism
between the category of (resp. nonlocal) vertex $F$-algebras and the
category of ordinary (resp. nonlocal) vertex algebras. We also
present some basic axiomatic properties similar to those for
ordinary vertex algebras.

We begin with the notion of nonlocal vertex algebra (see
\cite{li-qva1}; cf. \cite{bk}, \cite{li-g1}).

\bd{dnonlocalva} {\em A {\em nonlocal vertex algebra} is a vector
space $V$, equipped with a linear map
\begin{eqnarray*}
Y(\cdot,x):\ \ V&\rightarrow& \Hom (V,V((x)))\subset (\End V)[[x,x^{-1}]]\\
v&\mapsto& Y(v,x)=\sum_{n\in \Z}v_{n}x^{-n-1}\ \ (\mbox{where
}v_{n}\in \End V),
\end{eqnarray*}
and equipped with a distinguished vector ${\bf 1}\in V$, satisfying
the conditions that
\begin{eqnarray}\label{edef-vacuum}
 Y({\bf 1},x)v=v,\ \ Y(v,x){\bf
1}\in V[[x]]\ \ \mbox{ and }\ \lim_{x\rightarrow 0}Y(v,x){\bf 1}=v\
\ \mbox{ for }v\in V
\end{eqnarray}
(the {\em vacuum and creation property}), and that for $u,v,w\in V$,
there exists $l\in \N$ such that
\begin{eqnarray}\label{eweakassocrelation}
(x_{0}+x_{2})^{l}Y(u,x_{0}+x_{2})Y(v,x_{2})w
=(x_{0}+x_{2})^{l}Y(Y(u,x_{0})v,x_{2})w
\end{eqnarray}
(the {\em weak associativity}).} \ed

Note that in (\ref{eweakassocrelation}),  it is understood by
convention that
\begin{eqnarray*}
Y(u,x_{0}+x_{2})Y(v,x_{2})w&=&\sum_{m,n\in
\Z}(x_{0}+x_{2})^{-m-1}x_{2}^{-n-1}u_{m}v_{n}w\\
&=&\sum_{m,n\in \Z,\; j\ge
0}\binom{-m-1}{j}x_{0}^{-m-1-j}x_{2}^{-n-1+j}u_{m}v_{n}w,
\end{eqnarray*}
which exists in $V((x_{0}))((x_{2}))$. It is important to mention
that the expression
$$Y(u,x_{2}+x_{0})Y(v,x_{2})w$$ does {\em not} exist as a formal
series.

Now, let $F$ be a (one-dimensional) formal group over $\C$, which is
fixed throughout this section. We have $F(x,y)=F(y,x)$. Set $f=\log
F$.

Recall from \cite{li-phi-module} (cf. \cite{fhl}) the iota-map
$\iota_{x_{1},x_{2}}$, which was defined to be the unique algebra
embedding of the fraction field of $\C[[x_{1},x_{2}]]$ into
$\C((x_{1}))((x_{2}))$ such that
$$\iota_{x_{1},x_{2}}|_{\C[[x_{1},x_{2}]]}=1.$$
As $F(x,y)\in \C[[x,y]]\subset \C((x))[[y]]$ with $F(x,0)=x$ a unit
in $\C((x))$, $F(x,y)$, identified as $\iota_{x,y}F(x,y)$, is a unit
in the algebra $\C((x))[[y]]$. For the same reason, $F(x,y)$,
identified as $\iota_{y,x}F(x,y)$, is a unit in the algebra
$\C((y))[[x]]$. As a {\em convention}, for $m\in \Z$ we set
\begin{eqnarray*}
&&F(x,y)^{m}=(\iota_{x,y}F(x,y))^{m}\in \C((x))[[y]],\\
&&F(y,x)^{m}=(\iota_{y,x}F(x,y))^{m}\in \C((y))[[x]].
\end{eqnarray*}
Furthermore, for any $A(x_{1},x_{2})\in (\End W)[[x_{1}^{\pm
1},x_{2}^{\pm 1}]]$ with $W$ a vector space, we define
$A(F(x_{0},x_{2}),x_{2})$ and $A(F(x_{2},x_{0}),x_{2})$ accordingly.
Note that if $A(x_{1},x_{2})\in \Hom (W,W((x_{1}))((x_{2})))$,
$$A(F(x_{0},x_{2}),x_{2})\ \ \mbox{ exists in }\Hom (W,W((x_{0}))((x_{2}))).$$
If $A(x_{1},x_{2})\in \Hom (W,W((x_{1},x_{2})))$, in addition we
have
$$A(F(x_{2},x_{0}),x_{2})\ \ \mbox{ exists in }\Hom
(W,W((x_{2}))[[x_{0}]]).$$
For convenience, we set
\begin{eqnarray*}
&&A(x_{1},x_{2})|_{x_{1}=F(x_{0},x_{2})}=A(F(x_{0},x_{2}),x_{2}), \\
&&A(x_{1},x_{2})|_{x_{1}=F(x_{2},x_{0})}=A(F(x_{2},x_{0}),x_{2}).
\end{eqnarray*}
For $A(x_{1},x_{2})\in \Hom (W,W((x_{1}))((x_{2})))$, we have
\begin{eqnarray}\label{edouble-sub}
\left(A(x_{1},x_{2})|_{x_{1}=F(x_{0},x_{2})}\right)|_{x_{0}=f^{-1}(f(x_{1})-f(x_{2}))}
=A(x_{1},x_{2}),
\end{eqnarray}
where the substitution $x_{0}=f^{-1}(f(x_{1})-f(x_{2}))$ means
$$x_{0}^{m}=\iota_{x_{1},x_{2}}\left(f^{-1}(f(x_{1})-f(x_{2}))\right)^{m}
\ \ \ \mbox{ for all }m\in \Z.$$

\bd{dfva} {\em  We define a notion of {\em nonlocal vertex
$F$-algebra}, using all the axioms in Definition \ref{dnonlocalva}
except that weak associativity is replaced by the property that
 for $u,v,w\in V$, there exists $l\in \N$ such that
\begin{eqnarray}\label{eweak-assoc}
F(x_{0},x_{2})^{l}Y(u,F(x_{0},x_{2}))Y(v,x_{2})w=
F(x_{0},x_{2})^{l}Y(Y(u,x_{0})v,x_{2})w
\end{eqnarray}
(the {\em weak $F$-associativity}).} \ed

Furthermore, we formulate a notion of vertex $F$-algebra as follows:

\bd{dvfa} {\em A {\em vertex $F$-algebra} is a nonlocal vertex
$F$-algebra $V$ satisfying the {\em weak commutativity:} For $u,v\in
V$, there exists $k\in \N$ such that
\begin{eqnarray}\label{eweakcomm}
(x_{1}-x_{2})^{k}Y(u,x_{1})Y(v,x_{2})=(x_{1}-x_{2})^{k}Y(v,x_{2})Y(u,x_{1}).
\end{eqnarray}}
\ed

For convenience we formulate two technical lemmas.

\bl{lsimple-fact} Let $W$ be a vector space, let $a(x),b(x)\in \Hom
(W,W((x)))$, and let $g(x)\in x\C[[x]]$ with $g'(0)\ne 0$. Let $k$
be a nonnegative integer. Then
$$(g(x_{1})-g(x_{2}))^{k}a(x_{1})b(x_{2})\in \Hom (W,W((x_{1},x_{2})))$$
if and only if
$$(x_{1}-x_{2})^{k}a(x_{1})b(x_{2})\in \Hom (W,W((x_{1},x_{2}))).$$
On the other hand,
$$(g(x_{1})-g(x_{2}))^{k}a(x_{1})b(x_{2})=(g(x_{1})-g(x_{2}))^{k}b(x_{2})a(x_{1})$$
if and only if
$$(x_{1}-x_{2})^{k}a(x_{1})b(x_{2})=(x_{1}-x_{2})^{k}b(x_{2})a(x_{1}).$$
\el

\begin{proof} As $g(x)\in x\C[[x]]$ with $g'(0)\ne 0$, we have
$g(x_{1})-g(x_{2})=(x_{1}-x_{2})h(x_{1},x_{2})$ for some
$h(x_{1},x_{2})\in \C[[x_{1},x_{2}]]$ with $h(0,0)=g'(0)\ne 0$. Then
$$(g(x_{1})-g(x_{2}))^{n}=(x_{1}-x_{2})^{n}h(x_{1},x_{2})^{n}$$
for $n\in \N$. Note that as $h(0,0)\ne 0$, $h(x_{1},x_{2})$ is a
unit in the algebra $\C[[x_{1},x_{2}]]$, i.e., there exists
$q(x_{1},x_{2})\in \C[[x_{1},x_{2}]]$ such that
$h(x_{1},x_{2})q(x_{1},x_{2})=1$. It then follows immediately.
\end{proof}

The following can be proved similarly:

\bl{lsimple-fact2} Let $V$ and $W$ be vector spaces and let
\begin{eqnarray*}
Y(\cdot,x):& V\rightarrow& \Hom (V,V((x))), \\
Y_{W}(\cdot,x):& V\rightarrow &\Hom (W,W((x)))
\end{eqnarray*}
be linear maps. Let $u,v\in V$ and let $k\in \N$. Assume
\begin{eqnarray*}
(g(x_{1})-g(x_{2}))^{k}Y_{W}(u,x_{1})Y_{W}(v,x_{2})\in \Hom
(V,V((x_{1},x_{2})))
\end{eqnarray*}
and
\begin{eqnarray*}
&&(g(F(x_{0},x_{2}))-g(x_{2}))^{k}Y_{W}(Y(u,x_{0})v,x_{2})\nonumber\\
&=&\left[(g(x_{1})-g(x_{2}))^{k}Y_{W}(u,x_{1})Y_{W}(v,x_{2})\right]|_{x_{1}=F(x_{2},x_{0})}
\end{eqnarray*}
for some $g(x)\in x\C[[x]]$ with $g'(0)\ne 0$. Then the same
relations hold for every such $g(x)$. \el

The following is an explicit connection between (resp. nonlocal)
vertex $F$-algebras and ordinary (resp. nonlocal) vertex algebras:

\bp{pmain} Let $F$ be a formal group over $\C$, let $g(x)\in
x\C[[x]]$ with $g'(0)=1$, and let $F_{g}$ be the formal group
defined by
$$F_{g}(x,y)=g^{-1}(F(g(x),g(y)))$$
as in Proposition \ref{passociate-relation}. Let $V$ be a (resp.
nonlocal) vertex $F$-algebra. For $v\in V$, set
$$Y_{g}(v,x)=Y(v,g(x)).$$
Then $(V,Y_{g},{\bf 1})$ carries the structure of a (resp. nonlocal)
vertex $F_{g}$-algebra. \ep

\begin{proof} As $g(x)\in x\C[[x]]$ with $g'(0)=1$, for $v\in V$ we have
$$Y_{g}(v,x)=Y(v,g(x))\in \Hom(V,V((x))).$$
 We also have $Y_{g}({\bf 1},x)v=Y({\bf
1},g(x))v=v$, and
$$Y_{g}(v,x){\bf 1}=Y(v,g(x)){\bf 1}=\sum_{n\ge 0}g(x)^{n}v_{-n-1}{\bf 1},$$
which implies
$$Y_{g}(v,x){\bf 1}\in V[[x]]\ \ \mbox{ and }\ \ \lim_{x\rightarrow
0}Y_{g}(v,x){\bf 1}=v.$$ Furthermore, let $u,v,w\in V$. There exists
$l\in \N$ such that
\begin{eqnarray*}
F(z_{0},z_{2})^{l}Y(u,F(z_{0},z_{2}))Y(v,z_{2})w
=F(z_{0},z_{2})^{l}Y(Y(u,z_{0})v,z_{2})w.
\end{eqnarray*}
Write $g^{-1}(x)=xh(x)$ with $h(x)\in \C[[x]]$. Multiplying the both
sides of the equation above by $h(F(z_{0},z_{2}))^{l}$ which lies in
$\C[[z_{0},z_{2}]]$, we get
\begin{eqnarray*}
(g^{-1}(F(z_{0},z_{2})))^{l}Y(u,F(z_{0},z_{2}))Y(v,z_{2})w
=(g^{-1}(F(z_{0},z_{2})))^{l}Y(Y(u,z_{0})v,z_{2})w.
\end{eqnarray*}
Substituting $z_{0}=g(x_{0})$ and $z_{2}=g(x_{2})$ into the above
equation, we obtain
\begin{eqnarray*}
&&(g^{-1}(F(g(x_{0}),g(x_{2}))))^{l}Y(u,F(g(x_{0}),g(x_{2})))Y(v,g(x_{2}))w\\
&=&(g^{-1}(F(g(x_{0}),g(x_{2}))))^{l}Y(Y(u,g(x_{0}))v,g(x_{2}))w,
\end{eqnarray*}
which is
\begin{eqnarray*}
&&F_{g}(x_{0},x_{2})^{l}Y(u,g(F_{g}(x_{0},x_{2})))Y(v,g(x_{2}))w\\
&=&F_{g}(x_{0},x_{2})^{l}Y(Y(u,g(x_{0}))v,g(x_{2}))w.
\end{eqnarray*}
Namely,
\begin{eqnarray*}
F_{g}(x_{0},x_{2})^{l}Y_{g}(u,F_{g}(x_{0},x_{2}))Y_{g}(v,x_{2})w
=F_{g}(x_{0},x_{2})^{l}Y_{g}(Y_{g}(u,x_{0})v,x_{2})w.
\end{eqnarray*}
This proves that $(V,Y_{g},{\bf 1})$ carries the structure of a
nonlocal vertex $F_{g}$-algebra.

Now, assume that $V$ is a vertex $F$-algebra. Let $u,v\in V$. There
exists a nonnegative integer $k$ such that (\ref{eweakcomm}) holds.
In view of Lemma \ref{lsimple-fact} we have
$$(g^{-1}(x_{1})-g^{-1}(x_{2}))^{k}Y(u,x_{1})Y(v,x_{2})
=(g^{-1}(x_{1})-g^{-1}(x_{2}))^{k}Y(v,x_{2})Y(u,x_{1}).$$ With
substitution $x_{1}=g(z_{1})$, $x_{2}=g(z_{2})$, we get
$$(z_{1}-z_{2})^{k}Y_{g}(u,z_{1})Y_{g}(v,z_{2})
=(z_{1}-z_{2})^{k}Y_{g}(v,z_{2})Y_{g}(u,z_{1}).$$ Thus
$(V,Y_{g},{\bf 1})$ is a vertex $F_{g}$-algebra.
\end{proof}

Let $V$ be an ordinary (resp. nonlocal) vertex algebra. For $v\in
V$, set
$$Y_{F}(v,x)=Y(v,f(x)).$$
It follows from Proposition \ref{pmain} that $(V,Y_{F},{\bf 1})$
carries the structure of a (resp. nonlocal) vertex $F$-algebra.
Denote this (resp. nonlocal) vertex $F$-algebra by $V_{F}$.

To summarize, in terms of categories we have:

\bt{tmain} Let $F$ be a formal group over $\C$ with $\log F=f$.
 The map ${\mathcal{F}}$, defined by
${\mathcal{F}}(V)=V_{F}$ for every (resp. nonlocal) vertex algebra
$V$ and ${\mathcal{F}}(\theta)=\theta$ for every homomorphism
$\theta$ of (resp. nonlocal) vertex algebras, is an isomorphism from
the category of (resp. nonlocal) vertex algebras and the category of
(resp. nonlocal) vertex $F$-algebras.\et

\begin{proof} Let $\theta: U\rightarrow V$ be a homomorphism
of nonlocal vertex algebras. We have $\theta({\bf 1})={\bf 1}$ and
$$\theta(Y_{F}(u,x)u')=\theta(Y(u,f(x))u')=Y(\theta(u),f(x))\theta(u')
=Y_{F}(\theta(u),x)\theta(u')$$ for $u,u'\in U$. Thus $\theta$ is
also a homomorphism of nonlocal vertex $F$-algebras from $U_{F}$ to
$V_{F}$. It follows that ${\mathcal{F}}$ is a functor from the
category of (resp. nonlocal) vertex algebras to the category of
(resp. nonlocal) vertex $F$-algebras. On the other hand, by
Proposition \ref{pmain}, for each (resp. nonlocal) vertex
$F$-algebra $(K,Y,{\bf 1})$, $(K,Y(\cdot,f^{-1}(x)),{\bf 1})$ is an
ordinary (resp. nonlocal) vertex algebra. Then we have a functor
{}from the category of (resp. nonlocal) vertex $F$-algebras to  the
category of (resp. nonlocal) vertex algebras, sending each nonlocal
vertex $F$-algebra $(K,Y,{\bf 1})$ to $(K,Y(\cdot,f^{-1}(x)),{\bf
1})$. It follows immediately that ${\mathcal{F}}$ is an isomorphism.
\end{proof}

\bex{exFborcherds} {\em We here generalize Borcherds' construction
of (nonlocal) vertex algebras. Suppose that $F$ is a formal group
over $\C$ with $\log F=f$. Let $A$ be an associative algebra (over
$\C$) with identity and let $D$ be a derivation of $A$. Define
\begin{eqnarray}
Y_{F}(a,x)b=(e^{f(x)D}a)b\ \ \ \mbox{ for }a,b\in A.
\end{eqnarray}
It is known (see \cite{b-va}, \cite{b-Gva}) that for $F=F_{\rm a}$,
$(A,Y_{F_{\rm a}},1)$  carries the structure of a nonlocal vertex
algebra, which is a vertex algebra if and only if $A$ is
commutative.  Then by Proposition \ref{pmain},  for a general $F$,
$(A,Y_{F},1)$ carries the structure of a nonlocal vertex
$F$-algebra. Furthermore, $(A,Y_{F},1)$ is a vertex $F$-algebra if
and only if $A$ is commutative.} \eex

Next we present another version of weak $F$-associativity.

\bp{pF-wassoc} In the definition of a nonlocal vertex $F$-algebra,
weak $F$-associativity can be replaced by the property that for
$u,v\in V$, there exists $k\in \N$ such that
\begin{eqnarray}\label{epre-Fassoc}
(x_{1}-x_{2})^{k}Y(u,x_{1})Y(v,x_{2})\in \Hom (V,V((x_{1},x_{2}))),
\end{eqnarray}
and
\begin{eqnarray}\label{eFassociativity}
&&\left((x_{1}-x_{2})^{k}Y(u,x_{1})Y(v,x_{2})\right)|_{x_{1}=F(x_{2},x_{0})}\nonumber\\
&=&(F(x_{2},x_{0})-x_{2})^{k}Y(Y(u,x_{0})v,x_{2}).
\end{eqnarray}
Note that as
$\left(Y(u,x_{1})Y(v,x_{2})\right)|_{x_{1}=F(x_{2},x_{0})}$ does not
exist as a formal series, (\ref{epre-Fassoc}) is a precondition for
(\ref{eFassociativity}) to make sense. \ep

\begin{proof} First, assuming the very property we prove weak $F$-associativity.
Let $u,v,w\in V$. Let $k\in \N$ be such that (\ref{epre-Fassoc}) and
(\ref{eFassociativity}) hold. In view of Lemma \ref{lsimple-fact},
we have
\begin{eqnarray*}
(f(x_{1})-f(x_{2}))^{k}Y(u,x_{1})Y(v,x_{2}) \in
\Hom(V,V((x_{1},x_{2}))).
\end{eqnarray*}
Furthermore, there exists another nonnegative integer $l$ such that
\begin{eqnarray*}
x_{1}^{l}(f(x_{1})-f(x_{2}))^{k}Y(u,x_{1})Y(v,x_{2})w \in
V[[x_{1},x_{2}]][x_{2}^{-1}],
\end{eqnarray*}
involving only nonnegative powers of $x_{1}$. Then
\begin{eqnarray*}
&&\left(x_{1}^{l}(f(x_{1})-f(x_{2}))^{k}
Y(u,x_{1})Y(v,x_{2})w\right)|_{x_{1}=F(x_{2},x_{0})}\\
&=&\left(x_{1}^{l}(f(x_{1})-f(x_{2}))^{k}
Y(u,x_{1})Y(v,x_{2})w\right)|_{x_{1}=F(x_{0},x_{2})}\\
&=&F(x_{0},x_{2})^{l}(f(F(x_{0},x_{2}))-f(x_{2}))^{k}
Y(u,F(x_{0},x_{2}))Y(v,x_{2})w\\
&=&F(x_{0},x_{2})^{l}f(x_{0})^{k} Y(u,F(x_{0},x_{2}))Y(v,x_{2})w.
\end{eqnarray*}
On the other hand, writing
$f(x_{1})-f(x_{2})=(x_{1}-x_{2})q(x_{1},x_{2})$ with
$q(x_{1},x_{2})\in \C[[x_{1},x_{2}]]$, using (\ref{eFassociativity})
we have
\begin{eqnarray*}
&&\left(x_{1}^{l}(f(x_{1})-f(x_{2}))^{k}
Y(u,x_{1})Y(v,x_{2})w\right)|_{x_{1}=F(x_{2},x_{0})}\\
&=&\left(x_{1}^{l}q(x_{1},x_{2})^{k}(x_{1}-x_{2})^{k}
Y(u,x_{1})Y(v,x_{2})w\right)|_{x_{1}=F(x_{2},x_{0})}\\
 &=&F(x_{2},x_{0})^{l}f(x_{0})^{k}Y(Y(u,x_{0})v,x_{2})w.
\end{eqnarray*}
Consequently, we obtain
\begin{eqnarray}\label{ealmost}
&&F(x_{0},x_{2})^{l}f(x_{0})^{k}Y(u,F(x_{0},x_{2}))Y(v,x_{2})w\nonumber\\
&=&F(x_{0},x_{2})^{l}f(x_{0})^{k}Y(Y(u,x_{0})v,x_{2})w,
\end{eqnarray}
recalling that $F(x_{0},x_{2})=F(x_{2},x_{0})$. By canceling the
factor $f(x_{0})^{k}$ we obtain (\ref{eweak-assoc}).

Now, assume weak $F$-associativity holds. Let $u,v\in V$. There
exists $k\in \N$ such that $x_{0}^{k}Y(u,x_{0})v\in V[[x_{0}]]$,
which implies $f(x_{0})^{k}Y(u,x_{0})v\in V[[x_{0}]]$. Let $w\in V$
be arbitrarily fixed. There exists $l\in \N$ such that
(\ref{eweak-assoc}) holds. Then we have (\ref{ealmost}). We see that
the left-hand side of (\ref{ealmost}) involves only nonnegative
powers of $x_{0}$, so does the right-hand side. As
\begin{eqnarray*}
&&\left(x_{1}^{l}(f(x_{1})-f(x_{2}))^{k}
Y(u,x_{1})Y(v,x_{2})w\right)|_{x_{1}=F(x_{0},x_{2})}\\
&=&F(x_{0},x_{2})^{l}f(x_{0})^{k}Y(u,F(x_{0},x_{2}))Y(v,x_{2})w,
\end{eqnarray*}
we have
\begin{eqnarray*}
\left(x_{1}^{l}(f(x_{1})-f(x_{2}))^{k}
Y(u,x_{1})Y(v,x_{2})w\right)|_{x_{1}=F(x_{0},x_{2})} \in
V[[x_{0},x_{2}]][x_{2}^{-1}].
\end{eqnarray*}
By substitution $x_{0}=f^{-1}(f(x_{1})-f(x_{2}))$ (recall
(\ref{edouble-sub})), we get
$$x_{1}^{l}(f(x_{1})-f(x_{2}))^{k}
Y(u,x_{1})Y(v,x_{2})w\in V[[x_{1},x_{2}]][x_{2}^{-1}].$$ Thus
$$(f(x_{1})-f(x_{2}))^{k}
Y(u,x_{1})Y(v,x_{2})w\in V((x_{1},x_{2})).$$ As $k$ depends on $u$
and $v$, but not $w$, we have
$$(f(x_{1})-f(x_{2}))^{k}
Y(u,x_{1})Y(v,x_{2})\in \Hom (V,V((x_{1},x_{2}))).$$  Furthermore,
we have
\begin{eqnarray*}
&&\left(x_{1}^{l}(f(x_{1})-f(x_{2}))^{k}
Y(u,x_{1})Y(v,x_{2})w\right)|_{x_{1}=F(x_{0},x_{2})}\\
&=&\left(x_{1}^{l}(f(x_{1})-f(x_{2}))^{k}
Y(u,x_{1})Y(v,x_{2})w\right)|_{x_{1}=F(x_{2},x_{0})}\\
&=&F(x_{2},x_{0})^{l}\left[(f(x_{1})-f(x_{2}))^{k}
Y(u,x_{1})Y(v,x_{2})w\right]|_{x_{1}=F(x_{2},x_{0})}.
\end{eqnarray*}
Combining this with (\ref{ealmost}) we get
\begin{eqnarray*}
&&F(x_{0},x_{2})^{l}f(x_{0})^{k}Y(Y(u,x_{0})v,x_{2})w\nonumber\\
&=&F(x_{0},x_{2})^{l}\left[(f(x_{1})-f(x_{2}))^{k}
Y(u,x_{1})Y(v,x_{2})w\right]|_{x_{1}=F(x_{2},x_{0})}.
\end{eqnarray*}
Multiplying both sides by $\iota_{x_{2},x_{0}}F(x_{0},x_{2})^{-l}$,
we obtain
\begin{eqnarray*}
&&f(x_{0})^{k}Y(Y(u,x_{0})v,x_{2})w\nonumber\\
&=&\left((f(x_{1})-f(x_{2}))^{k}
Y(u,x_{1})Y(v,x_{2})w\right)|_{x_{1}=F(x_{2},x_{0})}.
\end{eqnarray*}
With $k$ independent of $w$, we have (\ref{eFassociativity}).
\end{proof}

\br{rnonlocal-va-wa2} {\em  Note that in the special case with
$F=F_{\rm a}$, Proposition \ref{pF-wassoc} (for ordinary nonlocal
vertex algebras) follows immediately from Lemma 2.9 of \cite{ltw}
with $W=V$ (the adjoint module) and $\sigma=1$ (the identity
automorphism).} \er

The following is an analog of a result in the theory of vertex
algebras (cf. \cite{li-g1}):

\bl{lD-operator} Let $V$ be a nonlocal vertex $F$-algebra. Define a
linear operator $\D$ on $V$ by
$$\D v=v_{-2}{\bf 1}=\left(\frac{d}{dx} Y(v,x){\bf 1}\right)|_{x=0}
\ \ \ \mbox{ for }v\in V.$$  Set $f=\log F$. Then
\begin{eqnarray}\label{eDbarcket-der}
&&Y(v,x){\bf 1}=e^{f(x)\D}v,\nonumber\\
&&[\D,Y(v,x)]=Y(\D v,x)=\frac{1}{f'(x)}\frac{d}{dx}Y(v,x)\ \ \
\mbox{ for }v\in V.
\end{eqnarray}
 \el

\begin{proof} For $v\in V$, set $\hat{Y}(v,x)=Y(v,f^{-1}(x))$.
By Proposition \ref{pmain}, $(V,\hat{Y},{\bf 1})$ is an ordinary
nonlocal vertex algebra. Let $\hat{\D}$ be the operator on $V$
defined by $\hat{D}v=\lim_{x\rightarrow
0}\frac{d}{dx}\hat{Y}(v,x){\bf 1}$ for $v\in V$. We have
$\hat{\D}=\D$ as
$$\hat{D}(v)=\lim_{x\rightarrow 0}\frac{d}{dx}\hat{Y}(v,x){\bf 1}
=\lim_{x\rightarrow 0}\frac{d}{dx}Y(v,f^{-1}(x)){\bf 1} =v_{-2}{\bf
1}=\D(v).$$ Using those $\D$-properties for nonlocal vertex algebras
(see \cite{li-g1}) we obtain
\begin{eqnarray*}
&&Y(v,x){\bf 1}=\hat{Y}(v,f(x)){\bf
1}=e^{f(x)\hat{\D}}v=e^{f(x)\D}v,\\
&&[\D, Y(v,x)]=[\hat{\D}, \hat{Y}(v,f(x))]=\hat{Y}(\D v,f(x))=Y(\D
v,x),
\end{eqnarray*}
and
\begin{eqnarray*}
Y(\D v,x)=\hat{Y}(\hat{\D}
v,f(x))=\left(\frac{d}{dz}\hat{Y}(v,z)\right)|_{z=f(x)}=\frac{1}{f'(x)}\frac{d}{dx}Y(v,x),
\end{eqnarray*}
as desired.
\end{proof}

The following is an analog of a well known result in vertex algebra
theory (see \cite{fhl}, \cite{dl}, \cite{li-local}):

\bp{pvfa-Ddef} A vertex $F$-algebra can be defined as a vector space
$V$ equipped with a linear map $Y(\cdot,x): V\rightarrow \Hom
(V,V((x)))$, a vector ${\bf 1}\in V$, and a linear operator $\D$ on
$V$, such that (\ref{edef-vacuum}), (\ref{eDbarcket-der}), and weak
commutativity  hold. \ep

\begin{proof} We only need to prove that
$V$ is a vertex $F$-algebra under the very assumptions. For $v\in
V$, set $\hat{Y}(v,x)=Y(v,f^{-1}(x))$, where $f(x)=\log F$. We have
$\hat{Y}({\bf 1},x)v=Y({\bf 1},f^{-1}(x))v=v$, and
$$\hat{Y}(v,x){\bf 1}=Y(v,f^{-1}(x)){\bf 1}\in V[[x]]
\ \ \mbox{ and }\ \ \lim_{x\rightarrow 0}\hat{Y}(v,x){\bf 1}=v,$$
noticing that $f^{-1}(x)\in x\C[[x]]$ with $(f^{-1})'(0)=1$. Let
$u,v\in V$. There exists $k\in \N$ such that (\ref{eweakcomm})
holds. {}From the second half of the proof of Proposition
\ref{pmain}, we have
\begin{eqnarray*}
(z_{1}-z_{2})^{k}\hat{Y}(u,z_{1})\hat{Y}(v,z_{2})
=(z_{1}-z_{2})^{k}\hat{Y}(v,z_{2})\hat{Y}(u,z_{1}).
\end{eqnarray*}
Furthermore, using (\ref{eDbarcket-der}) we get
\begin{eqnarray*}
&&[\D,\hat{Y}(v,x)]=[\D, Y(v,f^{-1}(x))]
=\left(\frac{1}{f'(z)}\frac{d}{dz}Y(v,z)\right)|_{z=f^{-1}(x)}\\
&=&\frac{d}{dx}Y(v,f^{-1}(x))=\frac{d}{dx}\hat{Y}(v,x).
\end{eqnarray*}
By Proposition 2.2.4 of \cite{li-local}, $(V,\hat{Y},{\bf 1})$ is a
vertex algebra. Now, it follows from Proposition \ref{pmain} that
$(V,Y,{\bf 1})$ is a vertex $F$-algebra.
\end{proof}

Next, we discuss a Jacobi identity for vertex $F$-algebras. Set
$f=\log F$ as before. Recall that $f(x)\in \C[[x]]$ with $f(0)=0$
and $f'(0)=1$. We have
\begin{eqnarray}
&&f(x_{0})^{-1}\delta\left(\frac{f(x_{1})-f(x_{2})}{f(x_{0})}\right)
-f(x_{0})^{-1}\delta\left(\frac{f(x_{2})-f(x_{1})}{-f(x_{0})}\right)
\nonumber\\
&&\hspace{2cm}
=f(x_{1})^{-1}\delta\left(\frac{f(x_{2})+f(x_{0})}{f(x_{1})}\right).
\end{eqnarray}
To see the existence in $\C[[x_{0}^{\pm 1},x_{1}^{\pm 1},x_{2}^{\pm
1}]]$ of the three terms, let use consider the first term. By
definition, we have
\begin{eqnarray*}
&&f(x_{0})^{-1}\delta\left(\frac{f(x_{1})-f(x_{2})}{f(x_{0})}\right)
=\sum_{n\in \Z}f(x_{0})^{-n-1}(f(x_{1})-f(x_{2}))^{n}\\
&&\hspace{2cm}=\sum_{n\in \Z}\sum_{j\ge
0}\binom{n}{j}(-1)^{j}f(x_{0})^{-n-1}f(x_{1})^{n-j}f(x_{2})^{j},
\end{eqnarray*}
where $$f(x_{2})^{j}\in x_{2}^{j}\C[[x_{2}]],\ \ \ f(x_{1})^{n-j}\in
x_{1}^{n-j}\C[[x_{1}]],\ \ \  f(x_{0})^{-n-1}\in
x_{0}^{-n-1}\C[[x_{0}]].$$ We see that the expression exists, by
first considering the coefficient of a fixed power of $x_{2}$, to
reduce to finitely many $j$. Furthermore, for each fixed $j$, we
have
\begin{eqnarray*}
&&\sum_{n\ge 0}\binom{n}{j}f(x_{0})^{-n-1}f(x_{1})^{n-j}\in
\C((x_{0}))((x_{1})),\\
 &&\sum_{n<0
}\binom{n}{j}f(x_{0})^{-n-1}f(x_{1})^{n-j}\in \C((x_{1}))((x_{0})).
\end{eqnarray*}

Using Proposition \ref{pmain} we immediately get:

\bp{pvfa} In the definition of a vertex $F$-algebra, weak
$F$-associativity and weak commutativity can be replaced by the
property that for $u,v\in V$,
\begin{eqnarray}
&&f(x_{0})^{-1}\delta\left(\frac{f(x_{1})-f(x_{2})}{f(x_{0})}\right)Y(u,x_{1})Y(v,x_{2})
\nonumber\\&&\hspace{2cm}
-f(x_{0})^{-1}\delta\left(\frac{f(x_{2})-f(x_{1})}{-f(x_{0})}\right)Y(v,x_{2})Y(u,x_{1})
\nonumber\\
&&
=f(x_{1})^{-1}\delta\left(\frac{f(x_{2})+f(x_{0})}{f(x_{1})}\right)Y(Y(u,x_{0})v,x_{2})
\end{eqnarray}
(the {\em Jacobi $F$-identity}), where $f=\log F$. \ep

\section{$\phi$-coordinated modules for a vertex
$F$-algebra}

In this section, we study $\phi$-coordinated (quasi) modules for a
vertex $F$-algebra with $\phi$ an associate of $F$. As the main
results of this section, we give a canonical connection between
$\phi$-coordinated (quasi) modules for nonlocal vertex algebras and
for nonlocal vertex $F$-algebras, and we give a general construction
of nonlocal vertex $F$-algebras and their $\phi$-coordinated (quasi)
modules.

We begin with the notion of module for a nonlocal vertex algebra
(see \cite{li-g1}).

\bd{drecall-module} {\em Let $V$ be a nonlocal vertex algebra. A
{\em $V$-module} is a vector space $W$ equipped with a linear map
\begin{eqnarray*}
Y_{W}(\cdot,x):\  V&\rightarrow& \Hom (W,W((x)))\subset (\End W)[[x,x^{-1}]]\\
 v&\mapsto&Y_{W}(v,x),
\end{eqnarray*}
satisfying the conditions that $$Y_{W}({\bf 1},x)=1_{W}\ (\mbox{the
identity operator on }W)$$ and that for $u,v\in V$, $w\in W$, there
exists $l\in \N$ such that
\begin{eqnarray}
(x_{0}+x_{2})^{l}Y_{W}(u,x_{0}+x_{2})Y_{W}(v,x_{2})w
=(x_{0}+x_{2})^{l}Y_{W}(Y(u,x_{0})v,x_{2})w.
\end{eqnarray}} \ed

We define a notion of {\em quasi $V$-module} (cf. \cite{li-qva1},
\cite{li-phi-module}) by replacing the above weak associativity with
the property that for $u,v\in V$, there exists a nonzero power
series $q(x_{1},x_{2})\in \C[[x_{1},x_{2}]]$ such that
$$q(x_{1},x_{2})Y_{W}(u,x_{1})Y_{W}(v,x_{2})\in \Hom
(W,W((x_{1},x_{2})))$$ and
$$\left(q(x_{1},x_{2})Y_{W}(u,x_{1})Y_{W}(v,x_{2})\right)|_{x_{1}=x_{2}+x_{0}}
=q(x_{0}+x_{2},x_{2})Y_{W}(Y(u,x_{0})v,x_{2}).$$

\br{rrecall-module} {\em {}From Lemma 2.9 of \cite{ltw}, in the
definition of a module for a nonlocal vertex algebra $V$, weak
associativity can be replaced by the property that for $u,v\in V$,
there exists $k\in \N$ such that
$$(x_{1}-x_{2})^{k}Y_{W}(u,x_{1})Y_{W}(v,x_{2})\in \Hom
(W,W((x_{1},x_{2})))$$ and
$$\left((x_{1}-x_{2})^{k}Y_{W}(u,x_{1})Y_{W}(v,x_{2})\right)|_{x_{1}=x_{2}+x_{0}}
=x_{0}^{k}Y_{W}(Y(u,x_{0})v,x_{2}).$$ In view of this, a $V$-module
is indeed a quasi $V$-module. } \er

Now, let $F$ be a (one-dimensional) formal group over $\C$, which is
fixed for the rest of this section. Set $f=\log F$.

\bd{dphi-module} {\em Let $V$ be a nonlocal vertex $F$-algebra and
let $\phi$ be an associate of $F$. We define a notion of {\em
$\phi$-coordinated quasi $V$-module}, using all the axioms in
Definition \ref{drecall-module} except that weak associativity is
replaced by the property that for $u,v\in V$, there exists
$q(x_{1},x_{2})\in \C[[x_{1},x_{2}]]$  such that
$q(\phi(x_{2},x_{0}),x_{2})\ne 0$,
$$q(x_{1},x_{2})Y_{W}(u,x_{1})Y_{W}(v,x_{2})\in
\Hom(W,W((x_{1},x_{2})))$$ and
\begin{eqnarray}
&&\left(q(x_{1},x_{2})Y_{W}(u,x_{1})Y_{W}(v,x_{2})\right)|_{x_{1}=\phi(x_{2},x_{0})}
\nonumber\\
&&=q(\phi(x_{2},x_{0}),x_{2})Y_{W}(Y(u,x_{0})v,x_{2})
\end{eqnarray}
(the {\em weak $\phi$-associativity}). We define a notion of {\em
$\phi$-coordinated $V$-module} by strengthening weak
$\phi$-associativity with $q(x_{1},x_{2})$ required to be a
polynomial of the form $(x_{1}-x_{2})^{k}$ with $k\in \N$.
Furthermore, we call an $F$-coordinated $V$-module a {\em
$V$-module}. } \ed

We have the following analog of Theorem \ref {tmain}:

\bp{pphi-module} Let $F$ be a formal group over $\C$ with $f=\log
F$, let $\phi$ be an associate of $F$, and let $g(x)\in x\C[[x]]$
with $g'(0)=1$. Set
$$\phi_{g}(x,z)=g^{-1}(\phi(g(x), g(z))),$$
an associate of $F_{g}$ (recall Proposition
\ref{passociate-relation}). Let $V$ be a nonlocal vertex $F$-algebra
and let $(W,Y_{W})$ be a $\phi$-coordinated (resp. quasi)
$V$-module. For $v\in V$, set
$$Y_{W}^{g}(v,x)=Y_{W}(v,g(x))\in (\End W)[[x,x^{-1}]].$$
Then $(W,Y_{W}^{g})$ carries the structure of a
$\phi_{g}$-coordinated (resp. quasi) $V_{g}$-module, where
 $V_{g}$ is a nonlocal vertex $F_{g}$-algebra
with $V_{g}=V$ as a vector space and with $Y_{g}(v,x)=Y(v,g(x))$ for
$v\in V$ (recall Proposition \ref{pmain}). \ep

\begin{proof} We shall consider the quasi case, while the non-quasi case
will be clear from the proof. We have
\begin{eqnarray*}
&&Y_{W}^{g}(v,x)=Y_{W}(v,g(x))\in \Hom (W,W((x)))\ \ \ \mbox{ for }v\in V,\\
&&Y_{W}^{g}({\bf 1},x)=Y_{W}({\bf 1},g(x))=1_{W}.
\end{eqnarray*}
For $u,v\in V$, there exists $q(z_{1},z_{2})\in \C[[z_{1},z_{2}]]$
such that $q(\phi(z_{2},z_{0}),z_{2})\ne 0$,
$$q(z_{1},z_{2})Y_{W}(u,z_{1})Y_{W}(v,z_{2})\in \Hom
(W,W((z_{1},z_{2}))),$$ and
$$\left(q(z_{1},z_{2})Y_{W}(u,z_{1})Y_{W}(v,z_{2})\right)|_{z_{1}=\phi(z_{2},z_{0})}
=q(\phi(z_{2},z_{0}),z_{2})Y_{W}(Y(u,z_{0})v,z_{2}).$$ With
substitution $z_{i}=g(x_{i})$ for $i=0,1,2$, we get
$$q(g(x_{1}),g(x_{2}))Y_{W}(u,g(x_{1}))Y_{W}(v,g(x_{2}))\in \Hom
(W,W((x_{1},x_{2})))$$ and
\begin{eqnarray*}
&&\left[q(g(x_{1}),g(x_{2}))
Y_{W}(u,g(x_{1}))Y_{W}(v,g(x_{2}))\right]|_{g(x_{1})=\phi(g(x_{2}),g(x_{0}))}\\
&=&q(\phi(g(x_{2}),g(x_{0})),g(x_{2}))Y_{W}(Y(u,g(x_{0}))v,g(x_{2})).
\end{eqnarray*}
Notice that the substitution $x_{1}=\phi_{g}(x_{2},x_{0})$ amounts
to $g(x_{1})=\phi(g(x_{2}),g(x_{0})).$ Then
\begin{eqnarray*}
&&\left[q(g(x_{1}),g(x_{2}))
Y^{g}_{W}(u,x_{1})Y^{g}_{W}(v,x_{2})\right]|_{x_{1}=\phi_{g}(x_{2},x_{0})}\\
&=&q(g(\phi_{g}(x_{2},x_{0})),g(x_{2}))Y_{W}^{g}(Y_{g}(u,x_{0})v,x_{2}).
\end{eqnarray*}
We have $q(g(x_{1}),g(x_{2}))\in \C[[x_{1},x_{2}]]$ such that
$$q(g(\phi_{g}(x,z)),g(x))=q(\phi(g(x),g(z)),g(x))\ne 0.$$
This proves that $(W,Y_{W}^{g})$ is a $\phi_{g}$-coordinated quasi
$V_{g}$-module.
\end{proof}

Recall that for a nonlocal vertex algebra $V$, $V_{F}$ is a nonlocal
vertex $F$-algebra with $V_{F}=V$ as a vector space and
$Y_{F}(v,x)=Y(v,f(x))$ for $v\in V$. With Proposition
\ref{pphi-module} the following is immediate:

\bc{cphi-module} Let $V$ be a nonlocal vertex algebra, let $F$ be a
formal group with $f=\log F$, and let $\phi$ be an associate of
$F_{\rm a}$. Set
$$\bar{\phi}(x,z)=f^{-1}(\phi(f(x), f(z))),$$
an associate of $F$. The map, sending $(W,Y_{W})$ to $(W,Y_{W}^{f})$
and sending each morphism $\theta$ to itself, is an isomorphism from
the category of $\phi$-coordinated (resp. quasi) $V$-modules to the
category of $\bar{\phi}$-coordinated (resp. quasi) $V_{F}$-modules.
\ec

The following is another connection between $\phi$-coordinated
modules for nonlocal vertex algebras and nonlocal vertex
$F$-algebras:

\bt{tsame-map}Let $F$ be a formal group over $\C$, let $\phi$ be an
associate of $F$, and let $g\in x\C[[x]]$ with $g'(0)=1$ . Set
$$\hat{\phi}_{g}(x,z)=\phi(x,g(z)),$$
an associate of $F_{g}$ by Proposition \ref{psimple-connection}. Let
$V$ be a nonlocal vertex algebra and let $W$ be a vector space
equipped with a linear map $Y_{W}(\cdot,x):V\rightarrow (\End
W)[[x,x^{-1}]]$. Then $(W,Y_{W})$ is a $\phi$-coordinated (resp.
quasi) $V$-module if and only if $(W,Y_{W})$ is a
$\hat{\phi}_{g}$-coordinated (resp. quasi) module for nonlocal
vertex $F_{g}$-algebra $(V,Y_{g},{\bf 1})$. \et

\begin{proof} We shall just consider the quasi case, as the non-quasi case
will be clear from the proof. Note that the $\phi$-coordinated quasi
$V$-module structure and the $\hat{\phi}_{g}$-coordinated quasi
$V_{g}$-module structure have the same requirement that
$Y_{W}(v,x)\in \Hom (W,W((x)))$ for $v\in V$ and $Y_{W}({\bf
1},x)=1_{W}$. Let $u,v\in V$. Assume
\begin{eqnarray*}
q(x_{1},x_{2})Y_{W}(u,x_{1})Y_{W}(v,x_{2}) \in \Hom
(W,W((x_{1},x_{2})))
\end{eqnarray*}
for some $q(x_{1},x_{2})\in \C[[x_{1},x_{2}]]$. We see that
$q(\phi(x_{2},x_{0}),x_{2})\ne 0$ if and only if
$q(\phi(x_{2},g(x_{0})),x_{2})\ne 0$. Furthermore, a
$\phi$-coordinated quasi $V$-module structure requires
\begin{eqnarray*}
\left(q(x_{1},x_{2})Y_{W}(u,x_{1})Y_{W}(v,x_{2})\right)_{x_{1}=\phi(x_{2},x_{0})}
=q(\phi(x_{2},x_{0}),x_{2})Y_{W}(Y(u,x_{0})v,x_{2}) ,
\end{eqnarray*}
while a $\phi_{g}$-coordinated quasi $V_{g}$-module structure
requires
\begin{eqnarray*}
\left(q(x_{1},x_{2})Y_{W}(u,x_{1})Y_{W}(v,x_{2})\right)_{x_{1}=\hat{\phi}_{g}(x_{2},x_{0})}
=q(\hat{\phi}_{g}(x_{2},x_{0}),x_{2})Y_{W}(Y_{g}(u,x_{0})v,x_{2}),
\end{eqnarray*}
which amounts to
\begin{eqnarray*}
&&\left(q(x_{1},x_{2})Y_{W}(u,x_{1})Y_{W}(v,x_{2})\right)|_{x_{1}=\phi(x_{2},g(x_{0}))}\\
&=&q(\phi(x_{2},g(x_{0})),x_{2})Y_{W}(Y(u,g(x_{0}))v,x_{2}).
\end{eqnarray*}
Then the equivalence between a $\phi$-coordinated quasi $V$-module
structure and a $\hat{\phi}_{g}$-coordinated quasi $V_{g}$-module
structure is immediate.
\end{proof}

Next, we give a construction of nonlocal vertex $F$-algebras and
their $\phi$-coordinated (quasi) modules, by using the results of
\cite{li-phi-module}. Let $W$ be a vector space over $\C$. Set
$$\E(W)=\Hom(W,W((x)))\subset (\End W)[[x,x^{-1}]].$$
A finite sequence $a^{1}(x),\dots,a^{r}(x)$ in $\E(W)$ is said to be
{\em quasi compatible} if there exists a nonzero power series
$p(x,y)\in \C[[x,y]]$ such that
\begin{eqnarray}
\left(\prod_{1\le i<j\le r}p(x_{i},x_{j})\right)a^{1}(x_{1})\cdots
a^{r}(x_{r})\in \Hom (W,W((x_{1},\dots,x_{r}))).
\end{eqnarray}
It is said to be {\em compatible} if the above containment relation
holds for $p(x,y)=(x-y)^{k}$ for some $k\in \N$.
 We say a subset $U$ of $\E(W)$ is {\em (resp. quasi) compatible} if any
finite sequence in $U$ is (resp. quasi) compatible.

Let $F$ be a formal group over $\C$ and let $\phi(x,z)$ be an
associate of $F$. We define a notion of {\em $\phi$-quasi
compatible} sequence (subset) by additionally requiring that
$p(\phi(x,z),x)\ne 0$ in the above definition. Notice that  if
$\phi(x,z)\ne x$, in view of Lemma \ref{lnonzero}, $\psi$-quasi
compatibility is simply the same as quasi compatibility.

\bd{dphi-action} {\em Let $a(x),b(x)\in \E(W)$ be such that
$(a(x),b(x))$ is $\phi$-quasi compatible, i.e., there exists
$p(x,y)\in \C[[x,y]]$ with $p(\phi(x,z),x)\ne 0$ such that
\begin{eqnarray}\label{epabx}
p(x_{1},x_{2})a(x_{1})b(x_{2})\in \Hom (W,W((x_{1},x_{2}))).
\end{eqnarray}
We define $a(x)_{n}^{\phi}b(x)\in \E(W)$ for $n\in \Z$ in terms of
the generating function
\begin{eqnarray*}
Y_{\E}^{\phi}(a(x),z)b(x)=\sum_{n\in \Z}a(x)_{n}^{\phi}b(x)z^{-n-1}
\end{eqnarray*}
 by
\begin{eqnarray}
Y_{\E}^{\phi}(a(x),z)b(x)
=p(\phi(x,z),x)^{-1}\iota_{x,z}\left(p(x_{1},x)a(x_{1})b(x)\right)|_{x_{1}=\phi(x,z)},
\end{eqnarray}
where $p(\phi(x,z),x)^{-1}$ denotes the inverse of $p(\phi(x,z),x)$
in $\C((x))((z))$.} \ed

The same argument in \cite{li-qva1} shows that
$Y_{\E}^{\phi}(a(x),z)b(x)$ is well defined; it is independent of
the choice of $p(x,y)$.

A $\psi$-quasi compatible subspace $U$ of $\E(W)$ is said to be {\em
$Y_{\E}^{\phi}$-closed} if
$$a(x)_{n}^{\phi}b(x)\in U\ \ \ \mbox{ for }a(x),b(x)\in U,\; n\in
\Z.$$

The following generalizes the corresponding result of
\cite{li-phi-module}:

\bt{tnew}  Let $F$ be a formal group over $\C$ and let $\phi(x,z)$
be an associate of $F$. Let $W$ be a vector space and let $U$ be a
$\phi$-quasi compatible subset of $\E(W)$.  Then there exists a
$Y_{\E}^{\phi}$-closed $\phi$-quasi compatible subspace which
contains $U$ and $1_{W}$. Denote by $\<U\>_{\phi}$ the smallest such
$Y_{\E}^{\phi}$-closed $\phi$-quasi compatible subspace. Then
$(\<U\>_{\phi},Y_{\E}^{\phi},1_{W})$ carries the structure of a
nonlocal vertex $F$-algebra and $W$ is a $\phi$-coordinated quasi
$\<U\>_{\phi}$-module with $Y_{W}(a(x),z)=a(z)$ for $a(x)\in
\<U\>_{\phi}$. \et

\begin{proof} Let $f$ be the logarithm of $F$. Set
$$\tilde{\phi}(x,z)=\phi(x,f^{-1}(z)).$$
By Proposition \ref{psimple-connection}, $\tilde{\phi}(x,z)$ is an
associate of $F_{\rm a}$. For any $q(x,y)\in \C[[x,y]]$, we see that
$q(\phi(x,z),x)\ne 0$ if and only if
$q(\tilde{\phi}(x,z),x)=q(\phi(x,f^{-1}(z)),x)\ne 0$. Then
$\phi$-quasi compatibility is the sam as $\tilde{\phi}$-quasi
compatibility and hence $U$ is $\tilde{\phi}$-quasi compatible. By
Theorem 4.11 of \cite{li-phi-module}, there exists a
$Y_{\E}^{\tilde{\phi}}$-closed $\tilde{\phi}$-quasi compatible
subspace containing $U$ and $1_{W}$, and
$(\<U\>_{\tilde{\phi}},Y_{\E}^{\tilde{\phi}},1_{W})$ carries the
structure of a nonlocal vertex algebra with $W$ as a
$\tilde{\phi}$-coordinated quasi module with $Y_{W}(a(x),z)=a(z)$
for $a(x)\in \<U\>_{\tilde{\phi}}$, where $\<U\>_{\tilde{\phi}}$ is
the smallest $Y_{\E}^{\tilde{\phi}}$-closed $\tilde{\phi}$-quasi
compatible subspace of $\E(W)$, containing $U$ and $1_{W}$.

Let $(a(x),b(x))$ be a $\phi$-quasi compatible pair in $\E(W)$.
There exists $p(x,y)\in \C[[x,y]]$ with $p(\phi(x,z),x)\ne 0$ such
that
$$p(x_{1},x_{2})a(x_{1})b(x_{2})\in \Hom (W,W((x_{1},x_{2}))).$$
We have
\begin{eqnarray*}
Y_{\E}^{\tilde{\phi}}(a(x),z)b(x)
&=&p(\tilde{\phi}(x,z),x)^{-1}\iota_{x,z}\left(p(x_{1},x)a(x_{1})b(x)\right)|_{x_{1}=\tilde{\phi}(x,z)},\\
Y_{\E}^{\phi}(a(x),z)b(x)
&=&p(\phi(x,z),x)^{-1}\iota_{x,z}\left(p(x_{1},x)a(x_{1})b(x)\right)|_{x_{1}=\phi(x,z)}.
\end{eqnarray*}
{}From this we get
\begin{eqnarray}\label{erelations}
Y_{\E}^{\phi}(a(x),z)b(x)=Y_{\E}^{\tilde{\phi}}(a(x),f(z))b(x).
\end{eqnarray}
By Theorem \ref{tmain}, $(\<U\>_{\tilde{\phi}},Y_{\E}^{\phi},1_{W})$
is a nonlocal vertex $F$-algebra. In particular,
$\<U\>_{\tilde{\phi}}$ is $Y_{\E}^{\phi}$-closed. This proves the
first assertion. It follows {}from (\ref{erelations}) that
$\<U\>_{\tilde{\phi}}$ is also the smallest $Y_{\E}^{\phi}$-closed
$\phi$-quasi compatible subspace containing $U$ and $1_{W}$. That
is, $\<U\>_{\phi}=\<U\>_{\tilde{\phi}}$.
 Furthermore, by Theorem
\ref{tsame-map}, $(W,Y_{W})$ is a $\tilde{\phi}$-coordinated quasi
$\<U\>_{\phi}$-module.
\end{proof}

\section{$\phi$-coordinated modules for nonlocal vertex algebras}

In this section, we focus our attention on $\phi$-coordinated
modules for nonlocal vertex algebras, in particular for ordinary
vertex algebras, with a specialized $\phi(x,z)=xe^{z}$. For a
nonlocal vertex $\Z$-graded algebra $V$, we exhibit a canonical
connection between $V$-modules and $\phi$-coordinated modules for a
nonlocal vertex algebra obtained from $V$ by Zhu's
change-of-variables theorem. Much of this section is motivated by
the work of Lepowsky (\cite{lep2}, \cite{lep3}).

Throughout this section, we fix
$$\phi(x,z)=xe^{z},$$
a particular associate of the additive formal group $F_{\rm a}$. We
have (recall Lemma \ref{lnonzero})
 $$q(xe^{z},x)\ne 0\ \  \mbox{ in } \C((x))[[z]]$$
for every nonzero $q(x,y)\in \C((x,y))$.

Let $V$ be a nonlocal vertex algebra. For a $\phi$-coordinated quasi
$V$-module $(W,Y_{W})$, the weak $\phi$-associativity states that
for any $u,v\in V$, there exists $0\ne q(x_{1},x_{2})\in
\C[[x_{1},x_{2}]]$ such that
\begin{eqnarray}
q(x_{1},x_{2})Y_{W}(u,x_{1})Y_{W}(v,x_{2})\in \Hom
(W,W((x_{1},x_{2})))
\end{eqnarray}
and
 \begin{eqnarray}
q(x_{2}e^{x_{0}},x_{2})Y_{W}(Y(u,x_{0})v,x_{2})
=\left(q(x_{1},x_{2})Y_{W}(u,x_{1})Y_{W}(v,x_{2})\right)|_{x_{1}=x_{2}e^{x_{0}}}.
\end{eqnarray}

The following was obtained in \cite{li-phi-module} (Lemma 3.2):

\bl{l-Dproperty-phi-module} Let $V$ be a nonlocal vertex algebra and
let $(W,Y_{W})$ be a $\phi$-coordinated quasi $V$-module. Then
\begin{eqnarray}
Y_{W}(e^{z\D} v,x)=Y_{W}(v,xe^{z})=e^{zx\frac{d}{dx}}Y_{W}(v,x) \ \
\ \mbox{ for }v\in V,
\end{eqnarray}
where $\D$ is the linear operator on $V$, defined by $\D
u=u_{-2}{\bf 1}$ for $u\in V$. In particular,
\begin{eqnarray}
Y_{W}(\D v,x)=x\frac{d}{dx}Y_{W}(v,x).
\end{eqnarray}
\el

While nonlocal vertex algebras are too general, what we called weak
quantum vertex algebras in \cite{li-qva1} form a special class of
nonlocal vertex algebras, which naturally generalize vertex algebras
and vertex superalgebras. A {\em weak quantum vertex algebra} is a
nonlocal vertex algebra $V$ satisfying the condition that for
$u,v\in V$, there exist
$$u^{(i)},\ v^{(i)}\in V,\ f_{i}(x)\in \C((x))\ \ (i=1,\dots,r)$$
such that
\begin{eqnarray}
&&x_{0}^{-1}\delta\left(\frac{x_{1}-x_{2}}{x_{0}}\right)
Y(u,x_{1})Y(v,x_{2})\nonumber\\
&&\ \ \ \ -x_{0}^{-1}\delta\left(\frac{x_{2}-x_{1}}{-x_{0}}\right)
\sum_{i=1}^{r}\iota_{x_{2},x_{1}}(f_{i}(x_{2}-x_{1}))
Y(v^{(i)},x_{2})Y(u^{(i)},x_{1})\nonumber\\
&=&x_{1}^{-1}\delta\left(\frac{x_{2}+x_{0}}{x_{1}}\right)
Y(Y(u,x_{0})v,x_{2}).
\end{eqnarray}

We recall the following result from \cite{li-phi-module}:

\bp{pwqva-phimodule} Let $V$ be a weak quantum vertex algebra and
let $(W,Y_{W})$ be a $\phi$-coordinated module for $V$ viewed as a
nonlocal vertex algebra. Let $u,v\in V$ and assume that
$$(x_{1}-x_{2})^{k}Y(u,x_{1})Y(v,x_{2})=(x_{1}-x_{2})^{k}
\sum_{i=1}^{r}\iota_{x_{2},x_{1}}(f_{i}(e^{x_{1}-x_{2}}))
Y(v^{(i)},x_{2})Y(u^{(i)},x_{1})$$ with $k\in \N,\; f_{i}(x)\in
\C(x),\; u^{(i)},v^{(i)}\in V$ for $1\le i\le r$. Then
\begin{eqnarray}\label{ejacobi-new4}
&&(x_{2}z)^{-1}\delta\left(\frac{x_{1}-x_{2}}{x_{2}z}\right)
Y_{W}(u,x_{1})Y_{W}(v,x_{2})\nonumber\\
&&\ \ \ \
-(x_{2}z)^{-1}\delta\left(\frac{x_{2}-x_{1}}{-x_{2}z}\right)
\sum_{i=1}^{r}\iota_{x_{2},x_{1}}(f_{i}(x_{1}/x_{2}))
Y_{W}(v^{(i)},x_{2})Y_{W}(u^{(i)},x_{1})\nonumber\\
&=&x_{1}^{-1}\delta\left(\frac{x_{2}(1+z)}{x_{1}}\right)
Y_{W}(Y(u,\log(1+z))v,x_{2}).
\end{eqnarray}
Furthermore, we have
\begin{eqnarray}\label{ejacobi-new5}
&&Y_{W}(u,x_{1})Y_{W}(v,x_{2})
-\sum_{i=1}^{r}\iota_{x_{2},x_{1}}(f_{i}(x_{1}/x_{2}))
Y_{W}(v^{(i)},x_{2})Y_{W}(u^{(i)},x_{1})\nonumber\\
&=&\Res_{x_{0}}\delta\left(\frac{x_{2}e^{x_{0}}}{x_{1}}\right)
Y_{W}(Y(u,x_{0})v,x_{2}).
\end{eqnarray}
 \ep

Note that
$$\delta\left(\frac{x_{2}e^{x_{0}}}{x_{1}}\right)
=e^{x_{0}\left(x_{2}\frac{\partial}{\partial
x_{2}}\right)}\delta\left(\frac{x_{2}}{x_{1}}\right).$$ As an
immediate consequence we have:

\bc{cphi-module-va} Let $V$ be a vertex algebra and let $(W,Y_{W})$
be a $\phi$-coordinated module for $V$ viewed as a nonlocal vertex
algebra. Then
\begin{eqnarray}
[Y_{W}(u,x_{1}),Y_{W}(v,x_{2})]=\sum_{j\ge
0}\frac{1}{j!}Y_{W}(u_{j}v,x_{2})\left(x_{2}\frac{\partial}{\partial
x_{2}}\right)^{j}\delta\left(\frac{x_{2}}{x_{1}}\right)
\end{eqnarray}
for $u,v\in V$. \ec

Furthermore, using Corollary \ref{cphi-module-va} and Lemma
\ref{l-Dproperty-phi-module} we immediately get:

\bc{cphi-module-voa} Let $(V,Y,{\bf 1},\omega)$ be a vertex operator
algebra of central charge $\ell\in \C$ in the sense of \cite{flm}
and let $(W,Y_{W})$ be a $\phi$-coordinated module for $V$ viewed as
a nonlocal vertex algebra. Then
\begin{eqnarray}
[L^{\phi}(m),L^{\phi}(n)]=(m-n)L^{\phi}(m+n)+\frac{\ell}{12}m^{3}\delta_{m+n,0}
\end{eqnarray}
for $m,n\in \Z$, where $Y_{W}(\omega,x)=\sum_{n\in
\Z}L^{\phi}(n)x^{-n}$. \ec

\bd{dzgradedva} {\em Let $\Gamma$ be a subgroup of the additive
group $\R$, containing $\Z$. A {\em nonlocal vertex $\Gamma$-graded
algebra} is a nonlocal vertex algebra $V$ equipped with a
$\Gamma$-grading $V=\bigoplus_{n\in \Gamma}V_{(n)}$, satisfying the
conditions that ${\bf 1}\in V_{(0)}$ and that for $v\in V_{(m)}$
with $m\in \Gamma$,
\begin{eqnarray}
v_{n}V_{(k)}\subset V_{(m+k-n-1)}\ \ \ \mbox{ for }n\in \Z,\; k\in
\Gamma.
\end{eqnarray}}
\ed

For a nonlocal vertex $\Gamma$-graded algebra $V$, define a linear
operator $L(0)$ on $V$ by
\begin{eqnarray}
L(0)|_{V_{(n)}}=n\ \ \ \mbox{ for }n\in \Gamma.
\end{eqnarray}
{}From \cite{fhl}, we have
\begin{eqnarray*}
L(0){\bf 1}&=&0,\\
x^{L(0)}Y(v,x_{1})x^{-L(0)}&=&Y(x^{L(0)}v,xx_{1}),\\
e^{xL(0)}Y(v,x_{1})e^{-xL(0)}&=&Y(e^{xL(0)}v,e^{x}x_{1})\ \ \ \mbox{
for }v\in V,
\end{eqnarray*}
where $x^{L(0)}u=x^{n}u$ for $u\in V_{(n)}$ with $n\in \Gamma$.

Furthermore, follow Zhu (\cite{zhu1}, \cite{zhu2}) to define a
linear map
$$Y[\cdot,x]: V\rightarrow (\End V)[[x,x^{-1}]]$$ by
\begin{eqnarray}
Y[v,x]=Y(e^{xL(0)}v,e^{x}-1)\ \ \ \mbox{ for }v\in V.
\end{eqnarray}

The following slightly generalizes Zhu's change-of-variables theorem
{}from vertex operator algebras to nonlocal vertex $\Gamma$-graded
algebras (cf. \cite{zhu1}, \cite{zhu2}):

\bp{pzhu-change} Let $V$ be a nonlocal vertex $\Gamma$-graded
algebra. Then $(V,Y[\cdot,x],{\bf 1})$ carries the structure of a
nonlocal vertex algebra. Furthermore, if $V$ is a vertex algebra,
$(V,Y[\cdot,x],{\bf 1})$ is also a vertex algebra. \ep

\begin{proof}  For $v\in V$, we have
$$Y[{\bf 1},x]v=Y(e^{xL(0)}{\bf 1},e^{x}-1)v=Y({\bf 1},e^{x}-1)v=v,$$
$$Y[v,x]{\bf 1}=Y(e^{xL(0)}v,e^{x}-1){\bf
1}=e^{(e^{x}-1)\D}e^{xL(0)}v\in V[[x]]$$ and $\lim_{x\rightarrow
0}Y[v,x]{\bf 1}=v$. Let $u,v,w\in V$. We have
\begin{eqnarray*}
Y[u,x_{0}+x_{2}]Y[v,x_{2}]w
&=&Y(e^{(x_{0}+x_{2})L(0)}u,e^{x_{0}+x_{2}}-1)Y(e^{x_{2}L(0)}v,e^{x_{2}}-1)w,\\
Y[Y[u,x_{0}]v,x_{2}]w
&=&Y(e^{x_{2}L(0)}Y(e^{x_{0}L(0)}u,e^{x_{0}}-1)v,e^{x_{2}}-1)w\\
&=&Y(Y(e^{(x_{2}+x_{0})L(0)}u,e^{x_{2}}(e^{x_{0}}-1))e^{x_{2}L(0)}v,e^{x_{2}}-1)w.
\end{eqnarray*}
Since $e^{(x_{0}+x_{2})L(0)}u\in V\otimes \C[[x_{0},x_{2}]]$ and
$e^{x_{2}L(0)}v\in V\otimes \C[[x_{2}]]$, there exists a nonnegative
integer $l$ such that
\begin{eqnarray*}
&&(z_{0}+z_{2})^{l}Y(e^{(x_{0}+x_{2})L(0)}u,z_{0}+z_{2})Y(e^{x_{2}L(0)}v,z_{2})w\\
&=&(z_{0}+z_{2})^{l}Y(Y(e^{(x_{2}+x_{0})L(0)}u,z_{0})e^{x_{2}L(0)}v,z_{2})w.
\end{eqnarray*}
Then by substituting $z_{0}=e^{x_{2}}(e^{x_{0}}-1), \;
z_{2}=e^{x_{2}}-1$ we obtain
$$(e^{x_{0}+x_{2}}-1)^{l}Y[u,x_{0}+x_{2}]Y[v,x_{2}]w
=(e^{x_{0}+x_{2}}-1)^{l}Y[Y[u,x_{0}]v,x_{2}]w,$$ noticing that
$(e^{x_{0}+x_{2}}-1)=e^{x_{2}}(e^{x_{0}}-1)+(e^{x_{2}}-1)$. As
$e^{x_{0}+x_{2}}-1=(x_{0}+x_{2})g(x_{0}+x_{2})$ for some $g(x)\in
\C[[x]]$ with $g(0)=1$, after cancellation we get
$$(x_{0}+x_{2})^{l}Y[u,x_{0}+x_{2}]Y[v,x_{2}]w
=(x_{0}+x_{2})^{l}Y[Y[u,x_{0}]v,x_{2}]w.$$ This proves that
$(V,Y[\cdot,x],{\bf 1})$ carries the structure of a nonlocal vertex
algebra.

Now, assume that $V$ is a vertex algebra. We shall prove that
$(V,Y[\cdot,x],{\bf 1})$ is also a vertex algebra by establishing
weak commutativity. Let $u,v\in V$. As $e^{xL(0)}u,\; e^{xL(0)}v\in
V\otimes \C[[x]]$, there exists $k\in \N$ such that
\begin{eqnarray*}
&&(y_{1}-y_{2})^{k}Y(e^{x_{1}L(0)}u,y_{1})Y(e^{x_{2}L(0)}v,y_{2})\\
&=&(y_{1}-y_{2})^{k}Y(e^{x_{2}L(0)}v,y_{2})Y(e^{x_{1}L(0)}u,y_{1}).
\end{eqnarray*}
Then
\begin{eqnarray*}
& &(e^{x_{1}}-e^{x_{2}})^{k}
Y(e^{x_{1}L(0)}u,e^{x_{1}}-1)Y(e^{x_{2}L(0)}v,e^{x_{2}}-1)\\
&=&(e^{x_{1}}-e^{x_{2}})^{k}Y(e^{x_{2}L(0)}v,e^{x_{2}}-1)
Y(e^{x_{1}L(0)}u,e^{x_{1}}-1).
\end{eqnarray*}
That is,
\begin{eqnarray*}
& &(e^{x_{1}}-e^{x_{2}})^{k}Y[u,x_{1}]Y[v,x_{2}]
=(e^{x_{1}}-e^{x_{2}})^{k}Y[v,x_{2}]Y[u,x_{1}].
\end{eqnarray*}
Writing $e^{x_{1}}-e^{x_{2}}=(x_{1}-x_{2})f(x_{1},x_{2})$, where
$f(x_{1},x_{2})=\sum_{n\ge
1}\frac{1}{n!}\frac{x_{1}^{n}-x_{2}^{n}}{x_{1}-x_{2}}$, we get
\begin{eqnarray*}
& &(x_{1}-x_{2})^{k}f(x_{1},x_{2})^{k} Y[u,x_{1}]Y[v,x_{2}]
=(x_{1}-x_{2})^{k}f(x_{1},x_{2})^{k}Y[v,x_{2}]Y[u,x_{1}].
\end{eqnarray*}
As $f(x_{1},x_{2})^{k}$ is invertible in $\C[[x_{1},x_{2}]]$, by
cancellation we obtain
\begin{eqnarray*}
(x_{1}-x_{2})^{k} Y[u,x_{1}]Y[v,x_{2}]
=(x_{1}-x_{2})^{k}Y[v,x_{2}]Y[u,x_{1}].
\end{eqnarray*}
Therefore, $(V,Y[\cdot,x],{\bf 1})$ is a vertex algebra.
\end{proof}

\br{risomorphism} {\em Let $(V,Y,{\bf 1},\omega)$ be a vertex
operator algebra of central charge $\ell\in \C$. Set
$\tilde{\omega}=\omega-\frac{\ell}{24}{\bf 1}$. It was proved in
\cite{zhu1} (cf. \cite{huang1}) that $(V,Y[\cdot,x],{\bf
1},\tilde{\omega})$ is a vertex operator algebra isomorphic to
$(V,Y,{\bf 1},\omega)$. } \er

The following is the main result of this section:

\bp{pordinary-exp-module} Let $V$ be a nonlocal vertex $\Z$-graded
algebra and let $(W,Y_{W})$ be a (resp. quasi) module for $V$ viewed
as a nonlocal vertex algebra. For $v\in V$, set
\begin{eqnarray}
X_{W}(v,x)=Y_{W}(x^{L(0)}v,x)\in \Hom (W,W((x)))\subset (\End
W)[[x,x^{-1}]].
\end{eqnarray}
Then $(W,X_{W})$ carries the structure of a $\phi$-coordinated
(resp. quasi) module with $\phi(x,z)=xe^{z}$ for the nonlocal vertex
algebra $(V,Y[\cdot,x], {\bf 1})$. \ep

\begin{proof} First, as $V$ is $\Z$-graded, we have
$X_{W}(v,x)w=Y_{W}(x^{L(0)}v,x)w\in W((x))$ for $v\in V,\; w\in W$.
Second, as $L(0){\bf 1}=0$ (with ${\bf 1}\in V_{(0)}$) we have
$$X_{W}({\bf 1},x)=Y_{W}(x^{L(0)}{\bf 1},x)=Y_{W}({\bf
1},x)=1_{W}.$$ Furthermore, let $u,v\in V$. As $x_{1}^{L(0)}u\in
V\otimes \C[x_{1},x_{1}^{-1}], \; x_{2}^{L(0)}v\in V\otimes
\C[x_{2},x_{2}^{-1}]$, it is clear that there exists a nonzero power
series $p(x_{1},x_{2})\in \C[[x_{1},x_{2}]]$ such that
$$p(x_{1},x_{2})X_{W}(u,x_{1})X_{W}(v,x_{2})\in \Hom
(W,W((x_{1},x_{2}))).$$ On the other hand, there exists a nonzero
power series $q(x_{1},x_{2})\in \C[[x_{1},x_{2}]]$ such that
$$q(x_{1},x_{2})Y_{W}((x_{2}e^{x_{0}})^{L(0)}u,x_{1})Y_{W}(x_{2}^{L(0)}v,x_{2})
\in \Hom (W,W((x_{1},x_{2})))[[x_{0}]],$$ and
\begin{eqnarray}\label{elast-one}
&& \left(q(x_{1},x_{2})Y_{W}((x_{2}e^{x_{0}})^{L(0)}u,x_{1})
Y_{W}(x_{2}^{L(0)}v,x_{2})\right)|_{x_{1}=x_{2}+z}\nonumber\\
&=&q(x_{2}+z,x_{2})Y_{W}(Y((x_{2}e^{x_{0}})^{L(0)}u,z)x_{2}^{L(0)}v,x_{2}).
\end{eqnarray}
We are going to apply substitution $z=x_{2}(e^{x_{0}}-1)$. Note that
for any $m\in \Z$,
\begin{eqnarray*}
&&\iota_{x_{2},x_{0}}\left[(x_{2}+z)^{m}|_{z=x_{2}(e^{x_{0}}-1)}\right]
=\sum_{i\ge 0}\binom{m}{i}x_{2}^{m-i}x_{2}^{i}(e^{x_{0}}-1)^{i}\\
&&\ \ \ \ =x_{2}^{m}\sum_{i\ge 0}\binom{m}{i}(e^{x_{0}}-1)^{i}
=x_{2}^{m}(e^{x_{0}})^{m}=x_{1}^{m}|_{x_{1}=x_{2}e^{x_{0}}}.
\end{eqnarray*}
Then
\begin{eqnarray*}
&&\left[ \left(q(x_{1},x_{2})Y_{W}((x_{2}e^{x_{0}})^{L(0)}u,x_{1})
Y_{W}(x_{2}^{L(0)}v,x_{2})\right)|_{x_{1}=x_{2}+z}\right]|_{z=x_{2}(e^{x_{0}}-1)}\\
&=&\left(q(x_{1},x_{2})Y_{W}(x_{1}^{L(0)}u,x_{1})
Y_{W}(x_{2}^{L(0)}v,x_{2})\right)|_{x_{1}=x_{2}e^{x_{0}}},
\end{eqnarray*}
while
\begin{eqnarray*}
&&\left[q(x_{2}+z,x_{2})
Y_{W}\left(Y((x_{2}e^{x_{0}})^{L(0)}u,z)x_{2}^{L(0)}v,x_{2}\right)\right]|_{z=x_{2}(e^{x_{0}}-1)}\\
&=&q(x_{2}e^{x_{0}},x_{2})Y_{W}\left(Y((x_{2}e^{x_{0}})^{L(0)}u,x_{2}(e^{x_{0}}-1))
x_{2}^{L(0)}v,x_{2}\right)\\
&=&q(x_{2}e^{x_{0}},x_{2})Y_{W}\left(x_{2}^{L(0)}Y(e^{x_{0}L(0)}u,e^{x_{0}}-1)
v,x_{2}\right)\\
&=&q(x_{2}e^{x_{0}},x_{2}) X_{W}(Y[u,x_{0}]v,x_{2}).
\end{eqnarray*}
Consequently,  by (\ref{elast-one}) we obtain
\begin{eqnarray*}
&&\left(q(x_{1},x_{2})X_{W}(u,x_{1})X_{W}(v,x_{2})\right)|_{x_{1}=x_{2}e^{x_{0}}}\\
&=&q(x_{2}e^{x_{0}},x_{2}) X_{W}(Y[u,x_{0}]v,x_{2}).
\end{eqnarray*}
This proves that $(W,X_{W})$ carries the structure of a
$\phi$-coordinated quasi module with $\phi(x,z)=xe^{z}$ for the
nonlocal vertex algebra $(V,Y[\cdot,x], {\bf 1})$. {}From the proof,
it is also clear for the non-quasi case.
\end{proof}

\br{rlast-comment} {\em Let $V$ be a vertex $\Z$-graded algebra and
let $(W,Y_{W})$ be a module for $V$ viewed as a nonlocal vertex
algebra. Combining Propositions \ref{pordinary-exp-module} and
\ref{pwqva-phimodule} we have
\begin{eqnarray}
[X_{W}(u,x_{1}),X_{W}(v,x_{2})]
=\Res_{x_{0}}\delta\left(\frac{x_{2}e^{x_{0}}}{x_{1}}\right)
X_{W}\left(Y[u,x_{0}]v,x_{2}\right)
\end{eqnarray}
for $u,v\in V$. This recovers a result of Lepowsky (\cite{lep3},
Theorem 3.1). In fact, it was the formal (unrigorous) relation
(2.18) in \cite{lep3} that motivated the study in this section.} \er

\end{document}